\documentclass[12pt,reqno]{amsart}
\usepackage{amsmath}
\usepackage{amssymb}
\usepackage{amstext}
\usepackage{mathrsfs}
\usepackage{a4wide}
\usepackage{graphicx}
\usepackage{bbm}
\usepackage{verbatim}
\usepackage{bm}
\usepackage{graphicx}
\usepackage{tikz}
\usepackage{pgfplots}

\allowdisplaybreaks \numberwithin{equation}{section}
\usepackage{color}
\usepackage{cases}
\usepackage{hyperref}
\hypersetup{hypertex=true,
	colorlinks=true,
	linkcolor=blue,
	anchorcolor=blue,
	citecolor=blue }

\numberwithin{equation}{section}

\newtheorem{theorem}{Theorem}[section]
\newtheorem{proposition}[theorem]{Proposition}

\newtheorem{lemma}[theorem]{Lemma}

\theoremstyle{definition}

\newtheorem{definition}[theorem]{Definition}

\theoremstyle{remark}
\newtheorem{remark}[theorem]{Remark}

\newcommand{\ep}{\varepsilon}

\begin{document}

\title
[local uniqueness of vortices for 2D steady Euler flow]{local uniqueness of vortices for 2D steady Euler flow in a bounded domain}

\author{Daomin Cao, Weilin Yu,  Changjun Zou}
\address{Institute of Applied Mathematics, Chinese Academy of Sciences, Beijing 100190, and University of Chinese Academy of Sciences, Beijing 100049,  P.R. China}
\email{dmcao@amt.ac.cn}

\address{Institute of Applied Mathematics, Chinese Academy of Sciences, Beijing 100190, and University of Chinese Academy of Sciences, Beijing 100049,  P.R. China}
\email{weilinyu@amss.ac.cn}

\address{Department of Mathematics, Sichuan University, Chengdu, Sichuan, 610064, P.R. China}
\email{zouchangjun17@mails.ucas.ac.cn}


\begin{abstract}
We study the 2D Euler equation in a bounded simply-connected domain, and establish the local uniqueness of flow whose stream function $\psi_\varepsilon$ satisfies
\begin{equation*}
	\begin{cases}
		-\varepsilon^2\Delta \psi_\varepsilon=\sum\limits_{i=1}^k \mathbf1_{B_\delta(z_{0,i})}(\psi_\varepsilon-\mu_{\varepsilon,i})_+^\gamma,\ \ \ &\text{in} \ \Omega,\\
		\psi_\varepsilon=0,\ \ \ & \text{on} \ \Omega, \\
	\end{cases}
\end{equation*}
with $\varepsilon\to 0^+$ the scale parameter of vortices, $\gamma\in(0,\infty)$, $\Omega\subset \mathbb R^2$ a bounded simply connected Lipschitz domain, $z_{0,i}\in\Omega$ the limiting location of $i^{\text{th}}$ vortex, and $\mu_{\varepsilon,i}$ the flux constants unprescribed. Our proof is achieved by a detailed description of asymptotic behavior for $\psi_\varepsilon$ and Pohozaev identity technique. For $k=1$, we prove the nonlinear stability of corresponding vorticity in $L^p$ norm, provided $z_{0,1}$ is a non-degenerate minimum point of Robin function. This stability result can be generalized to the case $k\ge 2$, and $(z_{0,1},\cdots,z_{0,k})\in \Omega^k$ being a non-degenerate minimum point of the Kirchhoff-Routh function.
\end{abstract}

\maketitle{\small{\bf Keywords:} The steady Euler equation; Kirchhoff-Routh function; local uniqueness; nonlinear stability. \\
	
}

\section{Introduction and main results}
In this paper, we study the planar flow of an ideal fluid, which is governed by the following Euler system
\begin{align}\label{1-1}
	\begin{cases}
		\partial_t\mathbf{u}+\mathbf{u}\cdot \nabla \mathbf{u} =-\nabla P,&\text{in}\ \Omega\times (0,T),\\
		\nabla\cdot\mathbf{u}=0,\,\ \, \ \ \ \ \ \ \  \ \, &\text{in}\  \Omega\times (0,T),\\
		\mathbf{u}\cdot \mathbf{n}=0, &\text{on}\  \partial \Omega,\\
		\mathbf{u}\big|_{t=0}=\mathbf{u}_0, &\text{in}\ \Omega,
	\end{cases}
\end{align}
where $\Omega\subset\mathbb{R}^2$ is a bounded simply connected domain with Lipschitz boundary, $\mathbf{u}=(u_1,u_2)$ is the velocity field, $P$ is the scalar pressure, $\mathbf{n}$ is the unit outward normal to $\partial \Omega$, and $\mathbf u_0$ is the initial data. 

The vorticity describing the rotation of the flow is defined by
\begin{equation*}
	\omega:=\nabla \times \mathbf{u}:=\partial_1u_2-\partial_2u_1.
\end{equation*}
By the incompressible condition $\nabla\cdot\mathbf u=0$ and slip boundary condition in \eqref{1-1}, there exists a Stokes stream function $\psi(x,t)$ such that 
$$\mathbf{u}=\nabla^\perp\psi:=(\partial_2\psi, -\partial_1\psi),$$ 
and $\psi=0$ on $\partial \Omega$. By the definition of $\psi$ and $\omega$, it holds $-\Delta \psi=\omega$. Hence we can recover the stream function $\psi$ by
\begin{equation*}
	\psi(x)=(-\Delta)^{-1}\omega(x)=\mathcal{G}\omega:=\int_{\Omega}G(x- y)\omega(y)dy,
\end{equation*}
where
$$G(x,y)=\frac{1}{2\pi}\ln \frac{1}{|x-y|}-H(x,y)$$ 
is the Green's function of $-\Delta$ in $\Omega$ with zero Dirichlet boundary condition, and $H(x,y)$ is the regular part. Now, we can rewrite \eqref{1-1} in the vorticity-stream function formulation
\begin{equation}\label{1-2}
	\begin{cases}
		\partial_t\omega+\mathbf{u}\cdot \nabla \omega =0,&\text{in}\ \Omega\times (0,T),
		\\
		\mathbf{u}=\nabla^\perp\psi,\,\,\, \psi =	(-\Delta)^{-1}\omega,&\text{in}\ \Omega\times (0,T),\\
		\omega\big|_{t=0}=\omega_0, &\text{in}\ \Omega,
	\end{cases}
\end{equation}
and it is easy to obtain the scalar pressure $P$ by $\mathbf u$ according to the first equation of \eqref{1-1}.

As a simplified version of 3D model, 2D incompressible Euler equation was intensively studied in the last century. In 1960s, Yudovich \cite{Yud} established the global well-posedness of \eqref{1-2} for $\omega$ in $L^1\cap L^\infty$. Since then, much effort has been made. DiPerna and Majda \cite{Dip} proved the existence of weak solutions for $\omega$ in $L^1\cap L^p$. In \cite{Del}, Delort investigated a general situation where the vorticity is a signed measure in $H^{-1}$. These results implies solutions to \eqref{1-2} is various. In this paper, we will first show the local uniqueness of some special steady solutions to \eqref{1-2}, where the vorticity $\omega$ is continuous with compact support. Then we discuss the nonlinear stability of the corresponding flow in some special situations.

To introduce our main approach, we define the weak solutions to \eqref{1-2} as follows.
\begin{definition}\label{def1}
	Given $\omega_0\in L^1(\Omega)\cap L^\infty(\Omega)$,  $\omega\in L^\infty([0,T]; L^1(\Omega)\cap L^\infty(\Omega))$ is called a weak solution to the Euler equation with the initial data $\omega_0$, provided that
	\begin{equation}\label{1-3}
		\int_{\Omega} \omega(x,T)\varphi(x,T)-\int_{\Omega} \omega_0(x)\varphi(x,0)=\int_0^T \int_{\Omega} \omega(\partial_t\varphi+\nabla^\perp 	\mathcal{G}\omega\cdot \nabla \varphi),\  \   \
	\end{equation}
	for all $\varphi\in C^1([0,T]; C_0^\infty(\Omega))$.
\end{definition}
For a steady flow we have $\omega\equiv \omega_0$, which is independent of the time parameter. Thus \eqref{1-3} is reduced to
\begin{equation}\label{1-4}
	\int_{\Omega} \omega\nabla^\perp \mathcal{G}\omega\cdot \nabla \varphi=0,\  \   \  \forall   \varphi\in C_0^\infty(\Omega).
\end{equation}
As mentioned in \cite{Bu2}, the general method of finding solutions to \eqref{1-4} is locally imposing $\omega=f(\psi)$ for some nondecreasing function $f:\mathbb R\to\mathbb R$. In our case, the stream function $\psi_\varepsilon$ satisfies
\begin{equation}\label{1-5}
	\begin{cases}
		-\varepsilon^2\Delta \psi_\varepsilon=\sum\limits_{i=1}^k \mathbf1_{B_\delta(z_{0,i})}(\psi_\varepsilon-\mu_{\varepsilon,i})_+^\gamma,\ \ \ &\text{in} \ \Omega,\\
		\psi_\varepsilon=0,\ \ \ & \text{on} \ \partial \Omega, \\
	\end{cases}
\end{equation}
where $\varepsilon>0$ is the scale parameter of vortices, $v_+=v$ if $v\ge0$, $v_+=0$ if $v<0$ for some function $v$,  $\gamma\in(0,\infty)$ is the power index, $z_{0,i}\in\Omega$ for $i=1,\cdots,k$ are $k$ different points describing the limiting location of each local vortex, $\delta$ is a small constant such that $$\{x\in B_\delta(z_{0,i}) \, | \, \psi_\varepsilon-\mu_{\varepsilon,i}>0\}\subset\subset B_\delta(z_{0,i})\subset\Omega \quad \text{and} \quad B_\delta(z_{0,i})\cap B_\delta(z_{0,j})=\emptyset, \quad \text{if} \ i\neq j.$$ 
For $i=1,\cdots,k$, $\mu_{\varepsilon,i}$ are $k$ flux constants unprescribed. Actually, $(\psi_\varepsilon,\mu_{\varepsilon,1},\cdots,\mu_{\varepsilon,k})$ constitutes a solution to the free boundary problem \eqref{1-5}, where $\psi_\varepsilon$ together with $\mu_{\varepsilon,i}$ determines the boundary of $i^{\text{th}}$ vortex, see the references \cite{Cao2,CWZ,SV}.

For further discussion, let $D_{\varepsilon,i}:=\{x\in B_\delta(z_{0,i}) \, | \, \psi_\varepsilon-\mu_{\varepsilon,i}>0\}$ be $k$ local vorticity sets, and let $\int_{B_\delta(x_{0,i})}(\psi_\varepsilon-\mu_{\varepsilon,i})_+^\gamma dx$ be the local circulation. To derive the local uniqueness, we assume additionally that the vorticity set $\cup_{i=1}^kD_{\varepsilon,i}$ of flow satisfies 
\begin{itemize}
	\item[($\mathbf{A}$)]  As $\varepsilon\to 0$, for $i=1,\cdots,k$, $\text{diam} \,D_{\varepsilon,i}\le R_0\varepsilon$ with $R_0$ a uniform constant, and $D_{\varepsilon,i}$ shrinks to $z_{0,i}$ in the sense that $ \sup_{x\in D_{\varepsilon,i}} |x-z_{0,i}|\to 0$;
\end{itemize}	
and the local circulation satisfies
\begin{itemize}
	\item[($\mathbf{B}$)]  For $i=1,\cdots,k$, it holds $\int_{B_\delta(x_{0,i})}(\psi_\varepsilon-\mu_{\varepsilon,i})_+^\gamma dx=\kappa_i$, where $\kappa_1,\cdots,\kappa_k$ are $k$ positive constants.
\end{itemize}

Notice that in \eqref{1-5}, we let the vorticity $\omega_\varepsilon$ be functionally dependent on $\psi_\varepsilon$ by a $\gamma$-power type $f$, which has been a favorite in previous studies \cite{Cao2,SV}. Moreover, our choice is especially typical since the change of $\gamma$ will lead to regularity difference on $\omega_\varepsilon$: if $\gamma\in(1,\infty)$, then the solutions to \eqref{1-5} give classical steady solutions to \eqref{1-2} in view of standard elliptic estimate. As $\gamma$ goes to $0$, solutions to \eqref{1-5} lose the regularity, and the case $\gamma=0$ in \cite{Cao4} is called the `vortex patch' since the vorticity $\omega$ is the characteristic function of $\cup_{i=1}^k D_{\varepsilon,i}$. Another interesting situation is $\gamma=\infty$, where the vorticity $\omega_\varepsilon=\varepsilon^2e^{\psi_\varepsilon}$ is in $C^\infty$ and supported on the whole domain $\Omega$ (The local uniqueness result in this paper does not contain this Liouville case). Readers can refer to \cite{DW} for relevant discussion. It should be mentioned that solutions in all of these cases possess a better regularity than point vortex solutions, where the $i^{\text{th}}$ local vortex is a Dirac measure with mass $\kappa_i$. Hence the procedure of solving \eqref{1-5} with ($\mathbf{A}$) and ($\mathbf{B}$) are also known as the regularization of point vortices in the study of Euler equations.

The topic of this paper originates from the study of vortex rings in the 3D axi-symmetric case \cite{AS,BF1,FT}. Thanks to the similarities of 3D axi-symmetric and planar Euler equations, Turkington \cite{T} constructed a family of 2D vortex patches, which approximate a single point vortex. Turkington's approach followed the lines of thought of Arnol'd \cite{Ar1,Ar2,Ar3}, and now is known as the vorticity method, whose strategy is to find the maximizers of kinetic energy in an $L^1\cap L^\infty$ admissible class of $\omega$. In view of the dual variational principle, these maximizers correspond to solutions to \eqref{1-4}. The second method for regularization is the stream function method, which focus on the semilinear elliptic equation that $\psi$ should satisfy (just like \eqref{1-5}), and deal with perturbed problem near ground states. Cao et al. adopted this method in \cite{Cao3}, and obtained multi-peak vortex patch solutions ($\gamma=0$)  to \eqref{1-4}. For the case $\gamma>0$, the study is carried out in a similar way. Readers can refer to \cite{CWZ} for vorticity method, and to \cite{Cao2,Cao1,SV} for stream function method. 

When $\gamma=1$, another source of problem \eqref{1-5} in reality is the plasma model (for detailed presentation of this model, see \cite{Tem}), which describes the equilibrium of a plasma confined in a toroidal cavity (the ``Tokamak machine"):
\begin{equation*}
	\begin{cases}
		\varepsilon^2\Delta u-u_-=0, \quad  & \text{in} \ \Omega,\\
		u=c_\Omega,\quad & \text{on} \ \partial\Omega, \\
		\int_{\partial\Omega}\frac{\partial u}{\partial\nu}=I,
	\end{cases}
\end{equation*}
where $\nu$ is the outward unit normal of $\partial\Omega$, $c_\Omega$ is a constant unprescribed, and $I$ is a given positive constant. With a solution $\psi_\varepsilon$ to \eqref{1-5} for $k=1$ and $\gamma=1$ in hand, we can use the linearity of equation and proceed as in \cite{Cao1} to obtain a solution for the plasma problem. Actually, if we let $\tilde u=1-\frac{\psi_\varepsilon}{\mu_{\varepsilon,1}}$, then for any $I$, $\tilde u_I=\frac{I}{\int_{\partial\Omega}\frac{\partial \tilde u}{\partial\nu}}\cdot\tilde u$ will be the desired solution. 

A long-standing question in history is the uniqueness of solutions to \eqref{1-5} as $\varepsilon\to0^+$, which can be obtained by two different ways introduced before: the variational method and the perturbed method. The significance of this question can be clarify as follows: once the local uniqueness is verified, one can follow the idea of Burton in \cite{Bu2}, and obtain the nonlinear stability of flow which maximizes the kinetic energy in a rearrangement class of vorticity $\omega_\varepsilon$, since the local uniqueness ensures the isolation of energy maximizers as a necessary condition for nonlinear stability. 
 
To begin with, let us determine the possible location for $z_{0,i}$ such that equation \eqref{1-5} is solvable. We consider the point vortex solution as a limiting case of $\varepsilon\to 0^+$, which is defined by
$$\omega_p(x)=\sum_{i=1}^{k}\kappa_i\boldsymbol\delta_{z_{0,i}}, \quad z_{0,i}\in\Omega$$
with $\boldsymbol\delta_x$ the Dirac measure centered at $x\in\Omega$. In \cite{Lin}, Lin showed that if $\omega_p$ is a steady solution to \eqref{1-2} in a very weak sense, then $\mathbf{z}_0=(z_{0,1},\cdots,z_{0,k})\in\Omega^k$ must be a critical point of the corresponding Kirchhoff-Routh function $\mathcal W_k$ for $\Omega$ and $k$ positive number $\kappa_1,\cdots,\kappa_k$, which is defined by
\begin{equation}\label{1-6}
	\mathcal{W}_k(x_1,x_2,..,x_k):=-\sum\limits_{i\neq j}^k\kappa_i\kappa_j G(x_i,x_j)+\sum\limits_{i=1}^k\kappa_i^2H(x_i,x_i).
\end{equation}
This is a direct consequence of $\mathcal W_k$ being the Hamiltonian dominating the evolution of $k$ point vortices. According to this fact, we may conjecture that the same condition should hold for $\varepsilon>0$ being very small, which is verified in the next theorem as the first result of this paper.

\begin{theorem}\label{thm1}
	Suppose for each $\varepsilon\in (0,\varepsilon_0]$ with $\varepsilon_0>0$ sufficiently small, \eqref{1-5} together with $(\mathbf{A})$ and $(\mathbf{B})$ has a solution $\psi_\varepsilon$. Then $\mathbf{z}_0=(z_{0,1},\cdots,z_{0,k})\in \Omega^k$ must be a critical point of $\mathcal{W}_k(\mathbf x)$.
\end{theorem}

The proof of Theorem \ref{thm1} is based on an asymptotic estimate for $\psi_\varepsilon$, and the following local Pohozaev identities for equation $-\Delta u=f(x,u)$ in $\Omega$:
\begin{equation}\label{1-7}
	\begin{split}
		&-\int_{\partial B_\delta(z_0)}\frac{\partial u}{\partial \nu}\frac{\partial u}{\partial x_h}+\frac{1}{2}\int_{\partial B_\delta(z_0)}|\nabla u|^2\nu_h\\
		&=\int_{\partial B_\delta(z_0)}F(x,u)\nu-\int_{B_\delta(z_0)}F_{x_i}(x,u), \ \ \ h=1,2,
	\end{split}
\end{equation}
with $x_0\in\Omega$, $\tau>0$ small, $\nu=(\nu_1,\cdots,\nu_N)$ the outward unit normal of $\partial B_\tau(x_0)$, and $F(x,t)=\int_0^{t}f(x,s)ds$. Actually, the Pohozaev identity technique is applied broadly in the study for `concentrating phenomena' of solutions to semilinear elliptic equations, for which we refer to Appendix 6.2 of book \cite{CPYb}.

Once Theorem \ref{thm1} is verified, we can further investigate the property of $\psi_\varepsilon$ by a more delicate estimate. To this purpose, for $\gamma\in (0,1)\cup(1,\infty)$, we denote $\phi_\gamma(y)=\phi_\gamma(|y|)$ as the unique radial solution of 
\begin{equation*}
	-\Delta\phi_\gamma=\phi_\gamma^\gamma, \ \ \phi_\gamma\in H^1_0(B_1(0)), \ \ \phi_\gamma>0 \ \ \text{in} \ B_1(0).
\end{equation*} 
While for $\gamma=1$, by Letting $\tau>0$ be the constant such that $1$ is the first eigenvalue of $-\Delta$ in $B_\tau(0)$ with the zero Dirichlet boundary condition, we denote $\phi_1(y)=\phi_1(|y|)>0$ as the first eigenfunction for $-\Delta$ in $B_\tau(0)$ with $\phi_1(0)=1$. We then have the following theorem on local uniqueness for $\psi_\varepsilon$ and some refined estimates. 
\begin{theorem}\label{thm2}
	Suppose that $\mathbf{z}_0=(z_{0,1},\cdots,z_{0,k})\in \Omega^k$ is a non-degenerate critical point of $\mathcal{W}_k(\mathbf{x})$. Then there exists a small number $\varepsilon_0>0$ such that for $\varepsilon\in (0,\varepsilon_0)$, equation \eqref{1-5} together with $(\mathbf{A})$ and $(\mathbf{B})$ has a unique solution $\psi_\varepsilon$.
	
	Moreover, for $i=1,\cdots,k$, it holds following asymptotic estimates:
	\begin{itemize}
		\item[(i)] The boundary of vorticity set $D_{\varepsilon,i}$ is a closed curve approximating a circle up to an error of order $O(\varepsilon^3)$, which can be parameterized as
		$$\partial D_{\varepsilon,i}=\left\{z_{\varepsilon,i}+r_{\varepsilon,i}(1+O(\varepsilon^2))(\cos(\theta),\sin(\theta))\,:\, \theta\in [0,2\pi)\right\},$$
		where
		\begin{equation*}
			z_{\varepsilon,i}=z_{0,i}+O(\varepsilon^2),
		\end{equation*}
	    \begin{equation*}
	    	r_{\varepsilon,i}=\varepsilon\cdot\left(\frac{2\pi\phi'_\gamma(1)}{\kappa_i}\right)^{\frac{\gamma-1}{2}}+O(\varepsilon^3), \quad \mathrm{for} \ \gamma\in(0,1)\cup(1,\infty),
	    \end{equation*}
        \begin{equation*}
        	r_{\varepsilon,i}=\tau\varepsilon+O(\varepsilon^3), \quad \mathrm{for} \ \gamma=1.
        \end{equation*}
		\item[(iii)] The flux constants
		\begin{equation*}
			\begin{split}
			\mu_{\varepsilon,i}&=\frac{\kappa_i}{2\pi} \ln \frac{1}{\varepsilon}+\frac{\gamma-1}{4\pi}\cdot\kappa_i \ln \frac{2\pi\phi_\gamma'(1)}{\kappa_i}-\kappa_iH(z_{0,i},z_{0,i})\\
			&\quad+\sum_{j\neq i}\kappa_j G(z_{0,i},z_{0,j})+O(\varepsilon^2|\ln\varepsilon|), \quad \mathrm{for} \ \gamma\in(0,1)\cup(1,\infty).
			\end{split}
		\end{equation*}
	    and
		\begin{equation*}
			\begin{split}
			\mu_{\varepsilon,i}&=\frac{\kappa_i}{2\pi} \ln \frac{1}{\varepsilon}+\frac{\kappa_i}{2\pi}\ln\tau-\kappa_iH(z_{0,i},z_{0,i})\\
			&\quad+\sum_{j\neq i}\kappa_j G(z_{0,i},z_{0,j})+O(\varepsilon^2|\ln\varepsilon|), \quad \mathrm{for} \ \gamma=1.
			\end{split}
		\end{equation*}
	\end{itemize}	
\end{theorem}

\begin{remark}
	In Section 3, we will show that the most important ingredient in the proof of local uniqueness is the $2$-dimensional kernel of the linearized operator. Indeed, for a general nonlinearity on the right hand side of the first equation in \eqref{1-5}, if we impose a non-degenerate condition on the corresponding ground state, then we can always expect a local uniqueness result just as Theorem \ref{thm2}. To achieve this goal, one can use a same argument as our proof for the case $\gamma=1$.
\end{remark}

We would also like to emphasize that in Theorem \ref{thm2}, $\mathbf{z}_0=(z_{0,1},\cdots,z_{0,k})\in \Omega^k$ is fixed as a non-degenerate critical point of $\mathcal{W}_k(\mathbf{x})$ and $\varepsilon\in(0,\varepsilon_0]$ is chosen to be very small. Hence our result is only concerning the local uniqueness of solutions. As a corollary of Theorem \ref{thm2}, if $\mathcal{W}_k(\mathbf{x})$ has more than one non-degenerate critical point, then there will be another different and locally unique $\psi_\varepsilon$ satisfying \eqref{1-5} by change our choice of $\mathbf{z}_0$.

As the scale parameter $\varepsilon$ becomes larger, the diameter of $D_{\varepsilon,i}$ will grow simultaneously. Since it is unclear whether the linearized operator has a fine property as the case $\varepsilon\in(0,\varepsilon_0]$, the global uniqueness of $\psi_\varepsilon$ remains to be a challenging open problem for general bounded domain $\Omega$. So far we only know the global uniqueness for special cases by symmetry, where $\Omega$ is an open disk and $k=1$, see \cite{Bu1} Theorem 3.1. On the other hand, if $\text{diam} \,D_{\varepsilon,i}$ exceeds $2\delta$ for some $1\le i\le k$, there will be no solution to \eqref{1-5}. Thus we can not expect the uniqueness of $\psi_\varepsilon$ for every $\varepsilon>0$. However, for a version of \eqref{1-5} where $k=1$ and the localization $\mathbf1_{B_\delta(z_{0,1})}$ is abandoned in the right hand side of equation, Bartolucci and Jevnikar \cite{Bar} proved that the solution $\psi_\varepsilon$ is unique for $\varepsilon>0$ very large (in this case the support of $\omega_\varepsilon$ will be the whole domain $\Omega$).  

To prove Theorem \ref{thm2}, we need a more precise estimate for $\psi_\varepsilon$ and vorticity set $\cup_{i=1}^k D_{\varepsilon,i}$, which is carried out by an approximation procedure and bootstrap. This approach can be also regarded as the inverse procedure of Lyapunov-Schmidt reduction for construction. Compared with the demonstration for patch case ($\gamma=0$) in \cite{Cao4}, where the authors calculate line integral along $\partial D_{\varepsilon,i}$, our estimates are simplified by the better regularity of solutions and an application of Taylor's formula. The argument in this paper also provides a possible approach for uniqueness of continuous vortex rings of small cross-section in 3D axi-symmetric flow.

Now we consider the special case $k=1$, where the Kirchhoff-Routh function $\mathcal W_1$ is also known as the Robin function
\begin{equation}\label{1-8}
	\mathcal R(x)=\kappa_1^2 H(x,x).
\end{equation}
To make the vorticity $\omega=-\Delta \psi$ maximize the kinetic energy of flow
$$E[\omega]=\frac{1}{2}\int_\Omega |\mathbf u|^2 dx=\frac{1}{2}\int_\Omega\omega\mathcal G\omega(x)dx$$
in a rearrangement class, we assume that the following condition holds:
\begin{itemize}
	\item[($\mathbf{C}$)] $k=1$, and $z_{0,1}$ is a non-degenerate minimum point of $\mathcal R(x)$.
\end{itemize}
According to Caffarelli and Friedman \cite{CF1}, if $\Omega$ is a convex domain, then $\mathcal R(x)$ is a strictly convex function, and $\mathcal R(x)$ has an isolated minimum in $\Omega$, which is also non-degenerate.

Using the local uniqueness result stated in Theorem \ref{thm2}, we can derive the nonlinear stability for a special solution $\omega_\varepsilon$ to \eqref{1-2} in $L^p$ norm, provided that $\psi_\varepsilon=(-\Delta)^{-1}\omega_\varepsilon$ is a solution to \eqref{1-5} together with conditions $(\mathbf{A})$, $(\mathbf{B})$ and $(\mathbf{C})$. In \cite{CWZ}, the authors showed that $\psi_\varepsilon$ will satisfy these assumptions if $\omega_\varepsilon=-\Delta\psi_\varepsilon$ is a maximizer of the kinetic energy relative to the rearrangement class generated by itself. Thus, in the spirit of Arnol'd \cite{Ar1,Ar2,Ar3}, Theorem \ref{thm2} will ensure that $\omega_\varepsilon$ is also a strict maximizer related to Kelvin's `isovortical surface’ \cite{K}, and hence it will constitute a stable steady flow (see Section 4 for discussion in detail).  

There are several known results on the Cauchy problem \eqref{1-2}: if $\omega\in L^{\infty}_{\mathrm{loc}}(\mathbb R, L^p(\Omega))$ with $4/3<p<\infty$ is a solution to \eqref{1-2}, then $\omega$ will also belong to the space $C(\mathbb R, L^p(\Omega))$, and the initial data $\omega_0$ can be defined. Moreover, for arbitrary $\omega_0\in L^p$ with $4/3<p<\infty$, although the uniqueness is absent by DiPerna and Majda \cite{Dip}, there will always exist global solutions that conserve kinetic energy. We refer to the appendix of Burton \cite{Bu2} for the rigorous 
proof of these statements. 

Having these preparations done, we are in the place to state the following theorem on nonlinear stability of vorticity.
\begin{theorem}\label{thm3}
	Let $\Omega\subset\mathbb{R}^2$ be a bounded simply connected domain of class $C^{2,\alpha}$ for some $\alpha\in(0,1)$, $\varepsilon\in (0,\varepsilon_0]$ with $\varepsilon_0>0$ sufficiently small, and $\psi_\varepsilon$ be solutions to \eqref{1-5} together with conditions $(\mathbf{A})$, $(\mathbf{B})$, $(\mathbf{C})$. Then $\omega_\varepsilon=-\Delta \psi_\varepsilon$ is nonlinearly stable in the following sense: for any $\eta>0$, there exists $\xi>0$ such that if $\omega\in L^{\infty}_{\mathrm{loc}}(\mathbb R, L^p(\Omega))$ with $4/3<p<\infty$ is an energy-conserving solution of \eqref{1-2}, and $\|\omega(0,\cdot)-\omega_\varepsilon\|_{L^p(\Omega)}<\xi$, then $\|\omega(t,\cdot)-\omega_\varepsilon\|_{L^p(\Omega)}<\eta$ for all $t\in\mathbb R$.
\end{theorem}
\begin{remark}
   If $\omega_0\in L^p$ with $3/2<p<\infty$, then all solutions to \eqref{1-2} conserve kinetic energy. Thus we do not have to assume $\omega$ is an energy-conserving solution in this situation. If $p=\infty$, Yudovich's theory \cite{Yud} guarantees the existence and uniqueness of solutions. However, Theorem \ref{thm3} does not apply to the $L^\infty$ norm.

   To easy our presentation, here we only consider the case $k=1$. However, for $k\ge2$ a finite natural number, $\mathbf{z}_0=(z_{0,1},\cdots,z_{0,k})\in \Omega^k$ a non-degenerate minimum point of $\mathcal{W}_k(\mathbf x)$, one can use a similar strategy to prove the nonlinear stability for $\omega_\varepsilon=-\Delta \psi_\varepsilon$ in $L^p$ norm. We will outline the proof in the end of Section 5.
\end{remark}

Using the uniqueness result obtained in \cite{Bu1}, Burton \cite{Bu2} gave an example of nonlinear stable vorticity where $\omega\in L^p(B_1(0))$ with $3/2<p<\infty$ is a non-negative radially symmetric decreasing function in an open disk. Theorem \ref{thm3} then can be regarded as a generalization of Burton's result to a bounded convex domain. For rotating and traveling-wave vortex patches to \eqref{1-2} in the whole $\mathbb R^2$, there are similar results on the nonlinear orbital stability, whose proof is also a combination of variational criterion and uniqueness argument, see \cite{Ab,CQZZ2}.  In addition, we would like to mention that Wan and Pulvirenti \cite{WP} proved the $L^1$ stability for circular vortex patch in a concentric circular domain. 
 
This paper is organized as follows. In Section 2, we study the asymptotic behavior of $\psi_\varepsilon$ and prove Theorem \ref{thm1}. In Section 3,  we improve the estimates by an approximation procedure and bootstrap. The uniqueness result stated in Theorem \ref{thm2} is verified in Section 4. In Section 5, we introduce the variational settings and prove Theorem \ref{thm3}.  Essential estimates for the free boundary $\partial D_{\varepsilon,i}$ are given in the Appendix.

\section{The necessary condition for the location of vortices}

In this section, we will prove Theorem \ref{thm1}, which relies on an accurate estimate of asymptotic behavior for $\psi_\varepsilon$. Let $\psi_\varepsilon$ be a solution to \eqref{1-5} together with conditions ($\mathbf{A}$) and ($\mathbf{B}$). By ($\mathbf{A}$), for $i=1,\cdots,k$ it holds 
$$\text{diam}\, D_{\varepsilon,i}\to 0, \quad \text{as} \ \varepsilon\to 0.$$ 
Let $r_{\varepsilon,i}=\frac{1}{2} \text{diam} \, D_{\varepsilon,i}$ be the radius of $D_{\varepsilon,i}$, and $p_{\varepsilon,i}\in D_{\varepsilon,i}$ be a point satisfying 
$$\psi_\varepsilon(p_{\varepsilon,i})=\max\limits_{x\in D_{\varepsilon,i}}\psi_\varepsilon(x),$$ 
which is always achievable by the maximum principle. We will start by giving estimates for the stream function $\psi_\varepsilon$ away from the vorticity set $\cup_{i=1}^k D_{\varepsilon,i}$. In the following, we will always assume that $L>0$ is a large constant independent of $\varepsilon$.
\begin{lemma}\label{lem2-1}
	Let $\psi_\varepsilon$ be a solution to \eqref{1-5} together with $(\mathbf{A})$ and $(\mathbf{B})$. For $x\in \Omega\setminus\cup^k_{i=1}\{x: \mathrm{dist}(x, D_{\varepsilon,i})\le Lr_{\varepsilon,i}\}$ it holds
	\begin{equation}\label{2-1}
	\psi_\varepsilon(x)=\sum\limits_{i=1}^k\kappa_i G(p_{\varepsilon,i},x)+O(\sum\limits_{i=1}^k\frac{ r_{\varepsilon,i}}{|x-p_{\varepsilon,i}|}),
	\end{equation}
	and 
	\begin{equation}\label{2-2}
	\frac{\partial\psi_\varepsilon(x)}{\partial x_h}=\sum\limits_{i=1}^k\kappa_i \frac{\partial G(p_{\varepsilon,i},x)}{\partial x_h}+O(\sum\limits_{i=1}^k\frac{ r_{\varepsilon,i}}{|x-p_{\varepsilon,i}|^2}), \quad h=1,2.
	\end{equation}
\end{lemma}
\begin{proof}
	For $x\in \Omega\setminus\cup^k_{i=1}\{x: \mathrm{dist}(x, D_{\varepsilon,i})\le Lr_{\varepsilon,i}\}$, it always holds $x \notin \cup_{i=1}^kD_{\varepsilon,i}$. According to the fact that
	$$H(y,x)-H(p_{\varepsilon,i},x)=O(\frac{r_{\varepsilon,j}}{|x-p_{\varepsilon,j}|}), \quad \forall \, y\in D_{\varepsilon,i},$$
	we deduce
	\begin{equation*}
	\begin{split}
	\psi_\varepsilon(x)&=\frac{1}{\varepsilon^2}\sum\limits_{i=1}^k\int_{D_{\varepsilon,i}}(\psi(y)-\mu_{\varepsilon,i})^\gamma_+G(y,x)dy\\
	&=\sum\limits_{i=1}^k\kappa_i G(p_{\varepsilon,i},x)+\frac{1}{\varepsilon^2}\sum\limits_{i=1}^k\int_{D_{\varepsilon,i}}(\psi(y)-\mu_{\varepsilon,i})^\gamma_+(G(y,x)-G(p_{\varepsilon,i},x))dy\\
	&=\sum\limits_{i=1}^k\kappa_i G(p_{\varepsilon,i},x)+O(\ln\frac{|x-p_{\varepsilon,i}|}{|y-x|})+O(\sum\limits_{i=1}^k\frac{ r_{\varepsilon,i}}{|x-p_{\varepsilon,i}|}).
	\end{split}
	\end{equation*} 
	Since
	\begin{equation*}
	|y-x|=|x-p_{\varepsilon, i}|-\langle\frac{x-p_{\varepsilon,i}}{|x-p_{\varepsilon,i}|},y-p_{\varepsilon,i}\rangle+O(\frac{|y-p_{\varepsilon,i}|^2}{|x-p_{\varepsilon,i}|}), \quad \forall \, y \in D_{\varepsilon,i},
	\end{equation*}
	we see \eqref{2-1} follows. Then the proof of \eqref{2-2} can be conducted in a similar way.
\end{proof}

To study the local behavior of $\psi_\varepsilon$ near $z_{\varepsilon,i}$, we will introduce the local limiting functions for $\psi_\varepsilon$ after scaling. For $\gamma\in (0,1)\cup(1,\infty)$, we let $\phi_\gamma(y)=\phi_\gamma(|y|)$ be the unique radial solution of 
\begin{equation*}
	-\Delta\phi_\gamma=\phi_\gamma^\gamma, \ \ \phi_\gamma\in H^1_0(B_1(0)), \ \ \phi_\gamma>0 \ \ \text{in} \ B_1(0).
\end{equation*}
Here, the uniqueness is obtained by an ODE technique, see \cite{GN,HS}. Then the local limiting function $w_\gamma$ is defined by
\begin{equation}\label{2-3}
	w_\gamma(y):=\left\{
	\begin{array}{lll}
		\phi_\gamma(|y|), & \mathrm{if} \ |y|\le 1,\\
		|\phi'_\gamma(1)|\ln\frac{1}{|y|}, & \mathrm{if} \ |y|>1.
	\end{array}
	\right.
\end{equation}
For $\gamma=1$, there will be a little bit different. We let $\tau>0$ be the constant such that $1$ is the first eigenvalue of $-\Delta$ in $B_\tau(0)$ with the zero Dirichlet boundary condition, and $\phi_1(y)=\phi_1(|y|)>0$ be the first eigenfunction for $-\Delta$ in $B_\tau(0)$ with $\phi_1(0)=1$. 
In this case, the local limiting function $w_1$ is given by
\begin{equation}\label{2-4}
	w_1(y):=\left\{
	\begin{array}{lll}
		\phi_1(|y|), & \mathrm{if} \ |y|\le \tau,\\
		\tau|\phi'_1(\tau)|\ln\frac{\tau}{|y|}, & \mathrm{if} \ |y|>\tau.
	\end{array}
	\right.
\end{equation}
Due to the importance of $w_\gamma$ in our argument, they are also known as the ``ground states" of the problem. Now we can state following proposition concerning the asymptotic estimates for $\psi_\varepsilon$.
\begin{proposition}\label{prop2-2}
	As $\varepsilon\to 0^+$, for $\gamma\in (0,1)\cup(1,\infty)$ and $i=1,\cdots,k$, one has
	\begin{equation*}
	\psi_\varepsilon(x)=\left(\frac{\varepsilon}{r_{\varepsilon,i}}\right)^{\frac{2}{\gamma-1}}\left(w\left(\frac{x-p_{\varepsilon,i}}{r_{\varepsilon,i}}\right)+o_\varepsilon(1)\right)+\mu_{\varepsilon,i}, \quad x\in B_{Lr_{\varepsilon,i}}(p_{\varepsilon,i}),
	\end{equation*}
    where $r_{\varepsilon,i}$, $\mu_{\varepsilon,i}$ satisfy
	\begin{equation}\label{2-5}
	\frac{\kappa_i}{2\pi}\ln\frac{1}{r_{\varepsilon,i}}+\sum\limits_{j\neq i}\kappa_jG(p_{\varepsilon,i},p_{\varepsilon,j})-\kappa_iH(p_{\varepsilon,i},p_{\varepsilon,j})-\mu_{\varepsilon,i}=o_\varepsilon(1),
	\end{equation}
	\begin{equation}\label{2-6}
	\kappa_i \left(\frac{r_{\varepsilon,i}}{\varepsilon}\right)^{\frac{2}{\gamma-1}}\to 2\pi|\phi'_\gamma(1)|.
	\end{equation}

	As $\varepsilon\to 0^+$, for $\gamma=1$ and $i=1,\cdots,k$, one has
	\begin{equation*}
	\psi_\varepsilon(x)=\Lambda \left(w\left(\frac{x-p_{\varepsilon,i}}{\varepsilon}\right)+o_\varepsilon(1)\right)+\mu_{\varepsilon,i}, \quad x\in B_{Lr_{\varepsilon,i}}(p_{\varepsilon,i}),
	\end{equation*}
    where $\Lambda$, $r_{\varepsilon,i}$, $\mu_{\varepsilon,i}$ satisfy
	\begin{equation}\label{2-7}
	\frac{\kappa_i}{2\pi}\ln\frac{1}{\tau\varepsilon}+\sum\limits_{j\neq i}\kappa_jG(p_{\varepsilon,i},p_{\varepsilon,j})-\kappa_iH(p_{\varepsilon,i},p_{\varepsilon,j})-\mu_{\varepsilon,i}=o_\varepsilon(1),
	\end{equation}
	\begin{equation}\label{2-8}
	\Lambda=\frac{\kappa_i}{2\pi\tau|\phi'_1(\tau)|} \quad \mathrm{and} \quad r_{\varepsilon,i}=\tau\varepsilon(1+o_\varepsilon(1)).
	\end{equation}
\end{proposition}

To prove Proposition \ref{prop2-2}, we need several lemmas to yield the convergence of scaled functions. In the first step, we obtain a uniform upper bound for the kinetic energy of fluid in vorticity set $\cup_{i=1}^kD_{\varepsilon,i}$.
\begin{lemma}\label{lem2-3}
	It holds
	$$\frac{1}{\varepsilon^2}\sum\limits_{i=1}^k\int_{D_{\varepsilon,i}}(\psi_\varepsilon-\mu_{\varepsilon,i})^{\gamma+1}_+dx= O_\varepsilon(1).$$
\end{lemma}
\begin{proof}
	According to local Pohozaev identities for \eqref{1-5} related to dilation, for $i=1,\cdots,k$, it holds
	\begin{equation*}
	\varepsilon^2\int_{\partial B_\delta(z_{0,i})}\langle x-p_{\varepsilon,i},\nabla\psi_\varepsilon\rangle\frac{\partial\psi_\varepsilon}{\partial\nu}-\frac{\varepsilon^2}{2}	\int_{\partial B_\delta(z_{0,i})}\langle x-p_{\varepsilon,j},\nu\rangle|\nabla\psi_\varepsilon|^2=\frac{2}{\gamma+1}\int_{D_{\varepsilon,i}}(\psi_\varepsilon-\mu_{\varepsilon,i})^{\gamma+1}_+dx
	\end{equation*}
	with $\nu$ the outward unit normal of $\partial B_\tau(z_{\varepsilon,i})$. Using Lemma \ref{lem2-1}, the left hand side of above identity is of order $O(\varepsilon^2)$. Hence we have verified the uniform bound.
\end{proof}

Now we let
$$w_{\gamma,\varepsilon,i}(y)=\left(\frac{r_{\varepsilon,i}}{\varepsilon}\right)^{\frac{2}{\gamma-1}}(\psi_\varepsilon(r_{\varepsilon,i}y+p_{\varepsilon,i})-\mu_{\varepsilon,i}), \quad \text{for} \ \gamma\in (0,1)\cup(1,\infty),$$
and
$$w_{1,\varepsilon,i}(y)=\Lambda^{-1}(\psi_\varepsilon(\varepsilon y+p_{\varepsilon,i})-\mu_{\varepsilon,i}).$$
Then for $i=1,\cdots,k$ we have
\begin{equation}\label{2-9}
-\Delta w_{\gamma,\varepsilon,i}=1_{B_{\delta r_{\varepsilon,i}^{-1}}(0)}(w_{\gamma,\varepsilon,i})^\gamma_++f_{\gamma,\varepsilon,i}(y), \quad \text{in} \ \Omega_{\gamma,\varepsilon,i},
\end{equation}
with 
\begin{equation*}
	\Omega_{\gamma,\varepsilon,i}:=\left\{
	\begin{array}{lll}
		\{y:r_{\varepsilon,j}y+p_{\varepsilon,i}\in\Omega\},& \text{for} \ \gamma\in (0,1)\cup(1,\infty),\\
		\{y:\varepsilon y+p_{\varepsilon,i}\in\Omega\}, & \text{for} \ \gamma=1,
	\end{array}
	\right.
\end{equation*}
and the remainder
\begin{equation*}
f_{\gamma,\varepsilon,i}(y):=\left\{
\begin{array}{lll}
\sum\limits_{i\neq j}1_{B_{\delta r_{\varepsilon,i}^{-1}}(\frac{p_{\varepsilon,j}-p_{\varepsilon,i}}{r_{\varepsilon,i}})}(w_{\gamma,\varepsilon,i}-\left(\frac{r_{\varepsilon,i}}{\varepsilon}\right)^{\frac{2}{\gamma-1}}(\mu_{\varepsilon,j}-\mu_{\varepsilon,i}))^\gamma_+, & \text{for} \ \gamma\in (0,1)\cup(1,\infty),\\
\sum\limits_{i\neq j}1_{B_{\delta r_{\varepsilon,i}^{-1}}(\frac{p_{\varepsilon,j}-p_{\varepsilon,i}}{\varepsilon})}(w_{\gamma,\varepsilon,i}-\Lambda^{-1}(\mu_{\varepsilon,j}-\mu_{\varepsilon,i}))_+,  & \text{for}  \ \gamma=1.
\end{array}
\right.
\end{equation*}
Intuitively, as $\varepsilon\to 0^+$, we see that $w_{\gamma,\varepsilon,i}$ might converge to $w_\gamma$, which solves
\begin{equation*}
	-\Delta w_\gamma=(w_\gamma)^\gamma_+, \quad \text{in} \ \mathbb R^2.
\end{equation*}
To verify this, we will derive a uniform bound for the $L^\infty$ norm of $w_{\gamma,\varepsilon,i}$ as $\varepsilon\to 0^+$.

\begin{lemma}\label{lem2-4}
	For any $R>0$, there exists a constant $C_R>0$ depending on $R$, such that
	$$||w_{\gamma,\varepsilon,i}||_{L^\infty(B_R(0)\cap\Omega_{\gamma,\varepsilon,i})}\le C_R.$$
\end{lemma}
\begin{proof}
	We are first to prove
	\begin{equation*}
	\int_{B_{\delta r_{\varepsilon,i}^{-1}}(0)}(w_{\gamma,\varepsilon,i})^{\gamma+1}_+dx=O_\varepsilon(1).
	\end{equation*}
	According to Lemma \ref{lem2-3}, we have
	\begin{equation}\label{2-10}
	\frac{1}{\varepsilon^2}\int_{D_{\varepsilon,i}}(\psi_\varepsilon-\mu_{\varepsilon,i})^{\gamma+1}_+dx=\left(\frac{\varepsilon}{r_{\varepsilon,i}}\right)^{\frac{4}{\gamma-1}}\int_{B_{\delta r_{\varepsilon,i}^{-1}}(0)}(w_{\gamma,\varepsilon,i})^{\gamma+1}_+\le C.
	\end{equation}
	By condition $(\mathbf{B})$ in Section 1, it holds
	\begin{equation*}
	\kappa_i=\frac{1}{\varepsilon^2}\int_{D_{\varepsilon,i}}(\psi_\varepsilon-\mu_{\varepsilon,i})^\gamma_+dx\le\left(\frac{1}{\varepsilon^2}\int_{D_{\varepsilon,i}}(\psi_\varepsilon-\mu_{\varepsilon,i})^{\gamma+1}_+\right)^{\frac{\gamma}{\gamma+1}}\cdot\left(\frac{|D_{\varepsilon,i}|}{\varepsilon^2}\right)^{\frac{1}{\gamma+1}},
	\end{equation*}
	which yields 
	$$\kappa_i^{\gamma+1}\le C\cdot\frac{|D_{\varepsilon,i}|}{\varepsilon^2}\le C\cdot \left(\frac{r_{\varepsilon,i}}{\varepsilon}\right)^2.$$ 
	On the other hand, by condition $(\mathbf{A})$ it holds $2r_{\varepsilon,i}\le R_0\varepsilon$. Combining the estimates for $r_{\varepsilon,i}$ in both directions, from \eqref{2-10} we obtain 
	\begin{equation*}
		\int_{B_{\delta r_{\varepsilon,i}^{-1}}(0)}(w_{\gamma,\varepsilon,i})^{\gamma+1}_+dx=O_\varepsilon(1).
	\end{equation*}
    Then applying Moser iteration on \eqref{2-9} (Theorem 4.1 in \cite{Han}), we conclude
	\begin{equation}\label{2-11}
	||(w_{\gamma,\varepsilon,i})_+||_{L^\infty(B_R(0)\cap\Omega_{\gamma,\varepsilon,i})}\le C.
	\end{equation}

	Now, let $v_1$ be a solution of 
	\begin{equation*}
	\begin{cases}
	-\Delta v_1=1_{B_{\delta r_{\varepsilon,i}^{-1}}(0)}(w_{\gamma,\varepsilon,i})^\gamma_++f_{\gamma,\varepsilon,i}(y), &\text{in} \  B_R(0)\cap\Omega_{\gamma,\varepsilon,i},\\
	v_1=0,  &\text{on} \  \partial(B_R(0)\cap\Omega_{\gamma,\varepsilon,i}).\\
	\end{cases}
	\end{equation*}
	Then it holds $|v_1|\le C$. Thus the function $v_2:=w_{\gamma,\varepsilon,i}-v_1$ will satisfy $-\Delta v_2=0$ in $B_R(0)\cap\Omega_{\gamma,\varepsilon,i}$, and
	\begin{equation}\label{2-12}
	\sup\limits_{B_R(0)\cap\Omega_{\gamma,\varepsilon,i}}v_2\ge\sup\limits_{B_R(0)\cap\Omega_{\gamma,\varepsilon,i}}w_{\gamma,\varepsilon,i}-C\ge-C.
	\end{equation}
	Since $\sup\limits_{B_R(0)\cap\Omega_{\gamma,\varepsilon,i}}w_{\gamma,\varepsilon,i}\ge 0$, \eqref{2-11} leads to 
	$$\sup\limits_{B_R(0)\cap\Omega_{\gamma,\varepsilon,i}}v_2\le\sup\limits_{B_R(0)\cap\Omega_{\gamma,\varepsilon,i}}w_{\gamma,\varepsilon,i}+C\le M$$ for some sufficiently large $M>0$. Thus $M-v_2$ is a positive harmonic function. By the Harnack inequality, there exists a constant $L>0$, such that
	$$\sup\limits_{B_R(0)\cap\Omega_{\gamma,\varepsilon,i}}(M-v_2)\le L\inf\limits_{B_R(0)\cap\Omega_{\gamma,\varepsilon,i}}(M-v_2).$$
	Combining this relationship with \eqref{2-12}, we have
	$$\inf\limits_{B_R(0)\cap\Omega_{\gamma,\varepsilon,i}}v_2\ge M-LM+L\sup\limits_{B_R(0)\cap\Omega_{\gamma,\varepsilon,i}}v_2\ge-C.$$
	Hence we have completed the proof.
\end{proof}

\begin{remark}
	When $\gamma\in(0,1)$, the assumption $\text{diam}\, D_ {\varepsilon,i}\le R_0\varepsilon$ in condition $(\mathbf{A})$ is not necessary in the proof of Lemma \ref{lem2-4}, since we only need a positive lower bound for $r_{\varepsilon,i}/\varepsilon$ in this case.
\end{remark}

Now we can give the local limiting function $w_\gamma$ a first description.
\begin{lemma}\label{lem2-6}
	As $\varepsilon\to 0^+$, it holds $w_{\gamma,\varepsilon,i}\to w_{\gamma,i}^*$ in $C^1_{\mathrm{loc}}(\mathbb{R}^2)$, where $w_{\gamma,i}^*$ is a radial function such that $(w_{\gamma,i}^*(r))'<0$ for all $r>0$.
\end{lemma}
\begin{proof}
	For $\gamma\in(0,1)\cup(1,\infty)$ and $y\in(B_{\delta r_{\varepsilon,i}^{-1}}(0)\setminus B_L(0))\cap\Omega_{\gamma,\varepsilon,i}$, it follows from asymptotic estimates in Lemma \ref{lem2-1} that
	\begin{equation*}
	\begin{split}
	w_{\gamma,\varepsilon,i}(y)&=\left(\frac{r_{\varepsilon,i}}{\varepsilon}\right)^{\frac{2}{\gamma-1}}(\psi_\varepsilon(r_{\varepsilon,i}y+p_{\varepsilon,i})-\mu_{\varepsilon,i})\\
	&=\left(\frac{r_{\varepsilon,i}}{\varepsilon}\right)^{\frac{2}{\gamma-1}}\left(\sum\limits_{j=1}^k\kappa_jG(r_{\varepsilon,i}y+p_{\varepsilon,i},p_{\varepsilon,j})-\mu_{\varepsilon,i}+O\left(\frac{1}{L}\right)\right)\\
	&=\left(\frac{r_{\varepsilon,i}}{\varepsilon}\right)^{\frac{2}{\gamma-1}}\cdot\frac{\kappa_i}{2\pi}\ln\frac{1}{|y|}+\left(\frac{r_{\varepsilon,i}}{\varepsilon}\right)^{\frac{2}{\gamma-1}}\cdot\\
	& \ \ \ \left(\frac{\kappa_i}{2\pi}\ln\frac{1}{r_{\varepsilon,i}}+\sum\limits_{j\neq i}\kappa_jG(r_{\varepsilon,i}y+p_{\varepsilon,i},p_{\varepsilon,j})-\mu_{\varepsilon,i}-\kappa_iH(r_{\varepsilon,i}y+p_{\varepsilon,i},p_{\varepsilon,i})+O\left(\frac{1}{L}\right)\right).
	\end{split}
	\end{equation*}
	By the proof of Lemma \ref{lem2-4}, it holds $(r_{\varepsilon,i}/\varepsilon)^{\frac{2}{\gamma-1}}\le C$. Thus we may assume (up to a subsequence) that $(r_{\varepsilon,i}/\varepsilon)^{\frac{2}{\gamma-1}}\kappa_i\to t_{\gamma,i}\in [0,+\infty)$. Moreover, it holds $w_{\gamma,\varepsilon,i}(y)\le C$ for each $y\in B_R(0)\cap\Omega_{\gamma,\varepsilon,i}$, which implies
	\begin{equation*}
	\left(\frac{r_{\varepsilon,i}}{\varepsilon}\right)^{\frac{2}{\gamma-1}}\cdot\left(\frac{\kappa_i}{2\pi}\ln\frac{1}{r_{\varepsilon,i}}+\sum\limits_{j\neq i}\kappa_jG(p_{\varepsilon,i},p_{\varepsilon,j})-\mu_{\varepsilon,i}-\kappa_iH(p_{\varepsilon,i},p_{\varepsilon,i})\right)\to \alpha_{\gamma,i}\in [0,+\infty).
	\end{equation*}
	For $\gamma=1$, one can use a similar approach to derive
	\begin{equation*}
	\begin{split}
	w_{1,\varepsilon,i}(y)&=\Lambda^{-1}(\psi_\varepsilon(\varepsilon y+p_{\varepsilon,i})-\mu_{\varepsilon,i})\\
	&=\Lambda^{-1}\cdot\left(\sum\limits_{j=1}^k\kappa_jG(\varepsilon y+p_{\varepsilon,i},p_{\varepsilon,j})-\mu_{\varepsilon,i}+O\left(\frac{1}{L}\right)\right)\\
	&=\Lambda^{-1}\cdot\frac{\kappa_i}{2\pi}\ln\frac{1}{|y|}+\Lambda^{-1}\cdot\\
	& \ \ \left(\frac{\kappa_i}{2\pi}\ln\frac{1}{\varepsilon}+\sum\limits_{j\neq i}\kappa_jG(\varepsilon y+p_{\varepsilon,i},p_{\varepsilon,j})-\mu_{\varepsilon,i}-\kappa_iH(\varepsilon y+p_{\varepsilon,i},p_{\varepsilon,i})+O\left(\frac{1}{L}\right)\right),
	\end{split}
	\end{equation*}
	and
	\begin{equation*}
	\Lambda^{-1}\cdot\left(\frac{\kappa_i}{2\pi}\ln\frac{1}{\varepsilon}+\sum\limits_{j\neq i}\kappa_jG(p_{\varepsilon,i},p_{\varepsilon,j})-\mu_{\varepsilon,i}-\kappa_iH(p_{\varepsilon,i},p_{\varepsilon,i})\right)\to \alpha_{1,i}\in [0,+\infty).
	\end{equation*}
	
	By the definition of $p_{\varepsilon,1},\cdots,p_{\varepsilon,k}$ at the start of this section, as $\varepsilon\to 0^+$, they will tend to $z_{0,1},\cdots,z_{0,k}$, which are $k$ different inner points of $\Omega$. Then according to \eqref{2-9}, we have $w_{\gamma,\varepsilon,i}\to w_{\gamma,i}^*$ in $C^1_{\text{loc}}(\mathbb{R}^2)$, where $w_{\gamma,i}^*(y)$ satisfies
	\begin{equation}\label{2-13}
	\begin{cases}
    -\Delta w_{\gamma,i}^*=(w_{\gamma,i}^*)_+^\gamma, &\text{in} \ B_R(0),\\
    w_{\gamma,i}^*=\frac{t_{\gamma,i}}{2\pi}\ln\frac{1}{|x|}+\alpha_{\gamma,i}, &\text{in}  \ B_R(0)\setminus B_L(0),
	\end{cases}
	\end{equation}
	with $R\gg L \gg 1$ being two large constants. Since $-\Delta w_{\gamma,i}^*\ge 0$, $w_{\gamma,i}^*$ will attain its minimum at the boundary of $B_R(0)$. Thus it holds
	$$w_{\gamma,i}^*(y)\ge\frac{t_{\gamma,i}}{2\pi}\ln\frac{1}{R}+\alpha_{\gamma,i},\quad \forall \, y\in B_R(0).$$
	Noticing that $(w_{\gamma,i}^*)_+^\gamma$ is an increasing function on $w_{\gamma,i}^*$, we can use the method  of moving plane to conclude that the solutions of \eqref{2-13} must be radially symmetric, and $(w_\gamma^*(r))'<0$ for all $r>0$.
\end{proof}

Now we are prepared to derive the asymptotic behavior of $\psi_\varepsilon$.

{\bf Proof of Proposition \ref{prop2-2}:}
We first consider the case $\gamma\in(0,1)\cup(1,\infty)$. By the definition of $r_{\varepsilon,i}$, we can find a $y_{\varepsilon,i}$ such that $|y_{\varepsilon,i}|=1$ and $r_{\varepsilon,i}y_{\varepsilon,i}+p_{\varepsilon,i}\in \partial D_{\varepsilon,i}$. Thus it holds
$$w_{\gamma,i}^*(y)=w_\gamma(y)$$
with $w_\gamma(y)$ defined in \eqref{2-3}. Comparing this result with the second equation in \eqref{2-13}, we obtain $t_\gamma=2\pi|\phi'_\gamma(1)|$ and $\alpha_\gamma+O(1/L)=0$. Since $\alpha_\gamma$ does not dependent on $L$, and $O(1/L)\to 0$ as $L\to0$, it must hold $\alpha_\gamma=0$. Hence as $\varepsilon\to 0$, we have
\begin{equation*}
\frac{\kappa_i}{2\pi}\ln\frac{1}{r_{\varepsilon,i}}+\sum\limits_{j\neq i}\kappa_jG(p_{\varepsilon,i},p_{\varepsilon,j})-\kappa_iH(p_{\varepsilon,i},p_{\varepsilon,j})-\mu_{\varepsilon,i}=o_\varepsilon(1),
\end{equation*}
and
\begin{equation*}
\kappa_i \left(\frac{r_{\varepsilon,i}}{\varepsilon}\right)^{\frac{2}{\gamma-1}}\to 2\pi|\phi'_\gamma(1)|,
\end{equation*}
which are exactly \eqref{2-5} and \eqref{2-6}.

For the case $\gamma=1$, we can study the asymptotic behavior of $\psi_\varepsilon$ similarly, and obtain
$$w_{1,i}^*(y)=w_1(y)$$
with $w_1(y)$ defined in \eqref{2-4}. In view of \eqref{2-13}, we conclude $t_1=2\pi \tau|\phi'_1(\tau)|$ and $\alpha_1=\tau|\phi'_1(\tau)|\ln\tau$. As $\varepsilon\to0$, it holds
\begin{equation*}
\frac{\kappa_i}{2\pi}\ln\frac{1}{\tau\varepsilon}+\sum\limits_{j\neq i}\kappa_jG(p_{\varepsilon,i},p_{\varepsilon,j})-\kappa_iH(p_{\varepsilon,i},p_{\varepsilon,j})-\mu_{\varepsilon,i}=o_\varepsilon(1).
\end{equation*}
One also has
\begin{equation*}
\Lambda=\frac{\kappa_i}{2\pi\tau|\phi'_1(\tau)|} \quad \mathrm{and} \quad r_{\varepsilon,i}=\tau\varepsilon(1+o_\varepsilon(1)).
\end{equation*}
Thus we have verified \eqref{2-7} and \eqref{2-8}, and finished the proof of Proposition \ref{prop2-2}.\qed\\

Taking advantage of the local Pohozaev identity \eqref{1-7}, we can give the necessary condition on the location of vortices.

{\bf Proof of Theorem \ref{thm1}:} 
Applying the Pohozaev identity \eqref{1-7} on \eqref{1-5} for the $i^{\text{th}}$ local vortex, we have
\begin{equation*}
-\int_{\partial B_\delta(z_{0,i})}\frac{\partial\psi_\varepsilon}{\partial\nu}\frac{\partial\psi_\varepsilon}{\partial x_h}+\frac{1}{2}\int_{\partial B_\delta(z_{0,i})}|\nabla\psi_\varepsilon|^2\nu_h=0, \quad h=1,2,
\end{equation*}
where $\nu=(\nu_1,\nu_2)$ is the outward unit normal of $\partial B_\delta(p_{\varepsilon,i})$.
Using the estimates in Lemma \ref{lem2-1}, we obtain 
\begin{equation*}
\begin{split}
&-\sum\limits_{j=1}^k\sum\limits_{m=1}^k\int_{\partial B_\delta(z_{0,i})}\kappa_j\kappa_m\langle DG(p_{\varepsilon,m},x),\nu\rangle D_{x_h}G(p_{\varepsilon,j},x)\\
&+\frac{1}{2}\int_{\partial B_\delta(z_{0,i})}(|\sum\limits_{j=1}^k\kappa_jDG(p_{\varepsilon,m},x)|)^2\nu_h=O(\sum\limits_{i=1}^kr_{\varepsilon,i}), \quad h=1,2,
\end{split}
\end{equation*}
where $\delta$ can be chosen arbitrarily small as $\varepsilon\to0^+$ from $(\mathbf{A})$. In view of Appendix 6.2 in \cite{CPYb}, the above relationship is equivalent to 
\begin{equation*}
\nabla_{\mathbf x}\mathcal{W}_k(z_{\varepsilon,1},\cdots,z_{\varepsilon,k})=O(\sum\limits_{i=1}^kr_{\varepsilon,i}).
\end{equation*}
Hence the proof of Theorem \ref{thm1} is complete.\qed

\section{The refined estimates for solutions}

In the preceding section, we have derived the necessary condition on the location of vortices. However, it is not enough for our purpose of proving the local uniqueness for solutions to \eqref{1-5} with the asymptotic behavior stated in Proposition \ref{prop2-2}. To refine the estimates for $\psi_\varepsilon$, we reconsider problem \eqref{1-5}
\begin{equation*}
	\begin{cases}
		-\varepsilon^2\Delta \psi_\varepsilon=\sum\limits_{i=1}^k \mathbf1_{B_\delta(x_{0,i})}(\psi_\varepsilon-\mu_{\varepsilon,i})_+^\gamma,\ \ \ &x\in\Omega,\\
		\psi_\varepsilon=0,\ \ \ \ &x\in\partial \Omega, \\
	\end{cases}
\end{equation*}
together with conditions $(\mathbf{A})$ and $(\mathbf{B})$. Our strategy is to construct a sequence of approximate solutions, and calculate their difference with $\psi_\varepsilon$ as an error term $\omega_\varepsilon$. Then we will bootstrap, and improve the estimates for $\omega_\varepsilon$ so that the difference is within the desired margins of error.

For the case $\gamma\in(0,1)\cup(1,\infty)$, we consider the following problem in $\mathbb R^2$:
\begin{equation*}
\begin{cases}
-\varepsilon^2\Delta u=(u-\frac{a}{2\pi}\ln\frac{1}{\varepsilon})_+^\gamma, &\text{in} \ B_s(z),\\
u=\frac{a}{2\pi}\ln\frac{1}{\varepsilon}, &\text{on} \ \partial B_s(z),
\end{cases}
\end{equation*}
where $s$ is a radius parameter to be chosen later. It has a unique continuous solution $W_{\gamma,\varepsilon,z,a}(x)$, which can be written as
\begin{equation*}
W_{\gamma,\varepsilon,z,a}(x)=\left\{
\begin{array}{lll}
\frac{a}{2\pi}\ln\frac{1}{\varepsilon}+(\varepsilon / s)^{\frac{2}{\gamma-1}}\phi_\gamma(\frac{|x-z|}{s}), & |x-z|\le s,\\
\frac{a}{2\pi}\frac{|\ln\varepsilon|}{|\ln s|}\ln\frac{1}{|x-z|}, & |x-z|\ge s,
\end{array}
\right.
\end{equation*}
where $\phi_\gamma(x)=\phi(|x|)$ is defined before Proposition \ref{prop2-2} as the unique radial solution of 
\begin{equation*}
	-\Delta\phi_\gamma=\phi_\gamma^\gamma, \ \ \phi_\gamma\in H^1_0(B_1(0)), \ \ \phi_\gamma>0 \ \ \text{in} \ B_1(0).
\end{equation*}
To make $W_{\gamma,\varepsilon,z,a}(x)$ in $C^1(\mathbb R^2)$, we also impose
\begin{equation*}
\left(\frac{\varepsilon}{s}\right)^{\frac{2}{\gamma-1}}\phi'_\gamma(1)=\frac{a}{2\pi}\frac{|\ln\varepsilon|}{|\ln s|}.
\end{equation*}
Let $z_{\varepsilon,i}$ be located near $z_{0,i}$, and $a_{\varepsilon,i}$ be the height parameter to be determined. We want to approximate the stream function $\psi_\varepsilon$ by superposition of $W_{\gamma,\varepsilon,z_{\varepsilon,i},a_{\varepsilon,i}}$ with $i=1,\cdots,k$. However, since $W_{\gamma,\varepsilon,z_{\varepsilon,i},a_{\varepsilon,i}}(x)\neq 0$ on $\partial\Omega$, we need to make the projection
\begin{equation*}
U_{\gamma,\varepsilon,z_{\varepsilon,i},a_{\varepsilon,i}}(x)=W_{\gamma,\varepsilon,z_{\varepsilon,i},a_{\varepsilon,i}}(x)-a_{\varepsilon,i}\cdot\frac{|\ln\varepsilon|}{|\ln s_{\varepsilon,i}|}\cdot H(x,z_{\varepsilon,i}),
\end{equation*}
where $H(x,y)$ is the regular part of Green's function for $\Omega$.

Now we let
\begin{equation*}
\mathcal{U}_{\gamma,\varepsilon,\mathbf{z_\varepsilon,a_\varepsilon}}(x)=\sum\limits_{i=1}^k U_{\gamma,\varepsilon,z_{\varepsilon,i},a_{\varepsilon,i}}(x),
\end{equation*}
where for $i=1,\cdots,k$, $\mathbf z_\varepsilon:=(z_{\varepsilon,1},\cdots,z_{\varepsilon,k})$, $\mathbf a_\varepsilon:=(a_{\varepsilon,1},\cdots,a_{\varepsilon,k})$ and $\mathbf s_\varepsilon:=(s_{\varepsilon,1},\cdots,s_{\varepsilon,k})$ are chosen to solve following system
\begin{equation}\label{3-1}
	\begin{cases}
		\nabla\mathcal{U}_{\gamma,\varepsilon,\mathbf{z_\varepsilon,a_\varepsilon}}(p_{\varepsilon,i})=0,\\
		\frac{a_{\varepsilon,i}}{2\pi}\ln\frac{1}{\varepsilon}=\mu_{\varepsilon,i}+\frac{a_{\varepsilon,i}|\ln\varepsilon|}{|\ln s_{\varepsilon,i}|}\cdot H(z_{\varepsilon,i},z_{\varepsilon,i})-\sum\limits_{j\neq i}\frac{a_{\varepsilon,j}|\ln\varepsilon|}{|\ln s_{\varepsilon,j}|}\cdot G(z_{\varepsilon,i},z_{\varepsilon,j}),\\
		\left(\frac{\varepsilon}{s_{\varepsilon,i}}\right)^{\frac{2}{\gamma-1}}\phi'_\gamma(1)=\frac{a_{\varepsilon,i}}{2\pi}\frac{|\ln\varepsilon|}{|\ln s_{\varepsilon,i}|},
	\end{cases}
\end{equation}
which can be transformed into a fixed point problem $\mathcal H_{\gamma,\varepsilon}(\mathbf{z_\varepsilon,a_\varepsilon, s_\varepsilon})= (\mathbf{z_\varepsilon,a_\varepsilon, s_\varepsilon})$. Using condition $(\mathbf{A})$ and contraction mapping theorem, we can verify that \eqref{3-1} has a unique solution in a neighborhood of $(\mathbf z_{0}, \boldsymbol \kappa, \mathbf s_0)$, where $\boldsymbol \kappa=(2\pi\kappa_1,\cdots,2\pi\kappa_k)$ and $\mathbf s_0=(s_{0,1},\cdots,s_{0,k})$ with $s_{0,i}=(2\pi\phi'_\gamma(1)/\kappa_i)^{\frac{\gamma-1}{2}}\varepsilon$. By Proposition \ref{prop2-2} and direct calculation, we can further deduce that
\begin{equation}\label{3-2}
	\begin{cases}
	    |z_{\varepsilon,i}-p_{\varepsilon,i}|=O(\varepsilon^2),\\
		|\frac{a_{\varepsilon,i}}{2\pi}\ln\frac{1}{\varepsilon}-\mu_{\varepsilon,i}|=O_\varepsilon(1),\\
		|r_{\varepsilon,i}-s_{\varepsilon,i}|=o(\varepsilon).
	\end{cases}
\end{equation}

For the case $\gamma=1$, we consider a slightly different problem in $\mathbb R^2$:
\begin{equation*}
	\begin{cases}
		-\varepsilon^2\Delta u= (u-\frac{a}{2\pi}\ln\frac{1}{\varepsilon})_+, \ u>0, &\text{in} \ B_{\tau\varepsilon}(z),\\
		u=\frac{a}{2\pi}\ln\frac{1}{\varepsilon}, &\text{on} \ \partial B_{\tau\varepsilon}(z),
	\end{cases}
\end{equation*}
where $a$ a positive parameter, and $\tau$ is the constant such that $1$ is the first eigenvalue of $-\Delta$ in $B_\tau(0)$ with the zero Dirichlet boundary condition. This problem has a unique solution $W_{1,\varepsilon,z,a}(x)$ defined as
\begin{equation*}
	W_{1,\varepsilon,z,a}(x)=\left\{
	\begin{array}{lll}
		\frac{a}{2\pi}\frac{1}{\varepsilon}+\Lambda_0\cdot\phi_1(\frac{|x-z|}{\varepsilon}),  & |x-z|\le \tau\varepsilon,\\
		\frac{a}{2\pi}\frac{|\ln\varepsilon|}{|\ln\tau\varepsilon|}\ln\frac{1}{|x-z|}, & |x-z|\ge \tau\varepsilon,
	\end{array}
	\right.
\end{equation*}
where $\phi_1(y)=\phi_1(|y|)>0$ is the first eigenfunction for $-\Delta$ in $B_\tau(0)$ with $\phi_1(0)=1$. To make $W_{1,\varepsilon,z,a}(x)\in C^1(\mathbb R^2)$, we impose
\begin{equation*}
	\tau\Lambda_0\cdot\phi'_1(\tau)=\frac{a}{2\pi}\frac{|\ln\varepsilon|}{|\ln \tau\varepsilon|}. 
\end{equation*}
As we have done in the former case, we let 
\begin{equation*}
	U_{1,\varepsilon,z_{\varepsilon,i},a_{\varepsilon,i}}(x)=W_{1,\varepsilon,z_{\varepsilon,i},a_{\varepsilon,i}}-a_{\varepsilon,i}\cdot\frac{|\ln\varepsilon|}{|\ln \tau\varepsilon|}\cdot H(x,z_{\varepsilon,i}).
\end{equation*}
Then an approximation of $\psi_\varepsilon$ is given by
\begin{equation*}
	\mathcal{U}_{1,\varepsilon,\mathbf{z_\varepsilon,a_\varepsilon}}(x)=\sum\limits_{i=1}^k U_{1,\varepsilon,z_{\varepsilon,i},a_{\varepsilon,i}}(x),
\end{equation*}
where for $i=1,\cdots,k$, parameters $\mathbf z_\varepsilon$, $\mathbf a_\varepsilon$ and $\boldsymbol \Lambda_\varepsilon:=(\Lambda_{\varepsilon,1},\cdots,\Lambda_{\varepsilon,k})$ are chosen to satisfy
\begin{equation}\label{3-3}
	\begin{cases}
		\nabla\mathcal{U}_{1,\varepsilon,\mathbf{z_\varepsilon,a_\varepsilon}}(p_{\varepsilon,i})=0,\\
		\frac{a_{\varepsilon,i}}{2\pi}\ln\frac{1}{\varepsilon}=\mu_{\varepsilon,i}+\frac{a_{\varepsilon,i}|\ln\varepsilon|}{|\ln \tau\varepsilon|}\cdot H(z_{\varepsilon,i},z_{\varepsilon,i})-\sum\limits_{j\neq i}\frac{a_{\varepsilon,j}|\ln\varepsilon|}{|\ln\tau\varepsilon|}\cdot G(z_{\varepsilon,i},z_{\varepsilon,j}),\\
		\tau\Lambda_{\varepsilon,i}\cdot\phi'_1(\tau)=\frac{a_{\varepsilon,i}}{2\pi}\frac{|\ln\varepsilon|}{|\ln \tau\varepsilon|}.
	\end{cases}
\end{equation}
System \eqref{3-3} can be also uniquely solved as a fixed point problem $\mathcal H_{1,\varepsilon}(\mathbf{z_\varepsilon,a_\varepsilon}, \boldsymbol \Lambda_\varepsilon)= (\mathbf{z_\varepsilon,a_\varepsilon}, \boldsymbol\Lambda_\varepsilon)$ near $(\mathbf z_{0}, \boldsymbol \kappa, \boldsymbol\Lambda_0)$, where $\boldsymbol \kappa=(2\pi\kappa_1,\cdots,2\pi\kappa_k)$ and $\boldsymbol \Lambda_0=(\Lambda_{0,1},\cdots,\Lambda_{0,k})$. Here $\Lambda_{0,i}=\kappa_i/2\pi\tau\phi'_1(\tau)$. Moreover, we have the asymptotic estimates
\begin{equation}\label{3-4}
	\begin{cases}
		|z_{\varepsilon,i}-p_{\varepsilon,i}|=O(\varepsilon^2),\\
		|\frac{a_{\varepsilon,i}}{2\pi}\ln\frac{1}{\varepsilon}-\mu_{\varepsilon,i}|=O_\varepsilon(1),\\
		|\Lambda_{\varepsilon,i}-\frac{\kappa_i}{2\pi\tau\cdot\phi'_1(\tau)}|=o_\varepsilon(1)
	\end{cases}
\end{equation}
from Proposition \ref{prop3-2}.

Now we are going to estimate the error term of the approximation defined by
\begin{equation*}
\omega_\varepsilon:=\psi_\varepsilon-\mathcal{U}_{\gamma,\varepsilon,\mathbf{z_\varepsilon,a_\varepsilon}}
\end{equation*} 
\begin{lemma}\label{lem3-1}
	As $\varepsilon\to0^+$, it holds
	$$||\omega_\varepsilon||_{L^\infty(\Omega)}+\varepsilon||\nabla\omega_\varepsilon||_{L^\infty(\Omega)}=o_\varepsilon(1).$$
\end{lemma}
\begin{proof}
	For the case $\gamma\in(0,1)\cup(1,\infty)$, using Proposition \ref{prop2-2} and estimates \eqref{3-2}, we can derive
	\begin{equation*}
	\begin{split}
	&\quad||\omega_\varepsilon||_{L^\infty(\cup^k_{i=1}B_{Lr_{\varepsilon,i}(z_{\varepsilon,i})})}=||\psi_\varepsilon-\mathcal{U}_{\gamma,\varepsilon,\mathbf{z_\varepsilon,a_\varepsilon}}||_{L^\infty(\cup^k_{i=1}B_{Lr_{\varepsilon,i}(z_{\varepsilon,i})})}\\
	&=\sup_{\cup^k_{i=1}B_{Lr_{\varepsilon,i}(z_{\varepsilon,i})}}\left|\left(\frac{\varepsilon}{r_{\varepsilon,i}}\right)^{\frac{2}{\gamma-1}}\left(\phi_\gamma(\frac{|x-p_{\varepsilon,i}|}{r_{\varepsilon,i}})+o_\varepsilon(1)\right)-\left(\frac{\varepsilon}{s_{\varepsilon,i}}\right)^{\frac{2}{\gamma-1}}\phi_\gamma(\frac{|x-z_{\varepsilon,i}|}{s_{\varepsilon,i}})+o_\varepsilon(1)\right|\\
	&=o_\varepsilon(1).
	\end{split}
	\end{equation*}
	Notice that $H(x,z_{\varepsilon,i})$ is bounded from below in $B_\delta(z_{0,i})\setminus B_{Lr_{\varepsilon,i}}(z_{\varepsilon,i})$. According to Lemma \ref{lem2-1}, for $x\in B_\delta(z_{0,j})\setminus B_{Lr_{\varepsilon,i}}(z_{\varepsilon,i})$ we have
	\begin{equation*}
	\begin{split}
	\psi_\varepsilon(x)&\le\frac{\kappa_i}{2\pi}\ln\frac{1}{|x-z_{\varepsilon,i}|}+O_\varepsilon(1)\\
	&=\frac{\mu_{i,\varepsilon}}{|\ln\varepsilon|}\left(1+O(\frac{1}{|\ln\varepsilon|})\right)\ln\frac{1}{|x-z_{\varepsilon,i}|}+O_\varepsilon(1)\\
	&\le\mu_{\varepsilon,i}\left(1-\frac{\ln L}{|\ln\varepsilon|}\right)\left(1+O(\frac{1}{|\ln\varepsilon|})\right)+O_\varepsilon(1)\\
	&<\mu_{\mu,i}.
	\end{split}
	\end{equation*}
	As a result, we obtain
	$$-\Delta\omega_\varepsilon=0, \ \ \ \text{in}  \ \Omega\setminus\cup^k_{i=1}B_{Lr_{\varepsilon,i}}(z_{\varepsilon,i}).$$
	By the maximum principle, we conclude
	$$||\omega_\varepsilon||_{L^\infty(\Omega\setminus\cup^k_{i=1}B_{Lr_{\varepsilon,i}}(z_{\varepsilon,i}))}\le||\omega_\varepsilon||_{L^\infty(\cup^k_{i=1}\partial B_{Lr_{\varepsilon,i}}(z_{\varepsilon,i}))}=o_\varepsilon(1).$$
	Next, we consider the $L^\infty$ norm of $\nabla\omega_\varepsilon$. Then for $x\in \Omega$, by Green's function and the estimate for $||\omega_\varepsilon||_{L^\infty(\Omega)}$, we have
	\begin{equation*}
		\begin{split}
			|\varepsilon\nabla \omega_\varepsilon(x)|&=\Bigg|\frac{1}{\varepsilon}\sum\limits_{i=1}^k\int_\Omega  \nabla_xG(x,y) \left(1_{B_\delta(z_{0,i})}(\psi_\varepsilon-\mu_{\varepsilon,i})^\gamma_+-(W_{\gamma,\varepsilon,z_{\varepsilon,i},a_{\varepsilon,i}}-\frac{a_{\varepsilon,i}}{2\pi}\ln\frac{1}{\varepsilon})_+^\gamma\right)dx\Bigg|\\
			&\le \frac{o_\varepsilon(1)}{\varepsilon}\sum\limits_{i=1}^k\int_{B_{Lr_{\varepsilon,i}}(z_{\varepsilon,i})}  \frac{1}{|x-y|} dx+o(\varepsilon)\\
			&=o_\ep(1).
		\end{split}
	\end{equation*}
	Hence we have verified the case $\gamma\in(0,1)\cup(1,\infty)$. For $\gamma=1$, we can use a similar argument to derive the desired estimate.
\end{proof}

Now, we will further investigate the error term $\omega_\varepsilon$ by method of linearization. Define the linear operator
\begin{equation*}
\mathbb{L}_{\gamma,\varepsilon}v:=-\varepsilon^2\Delta v-\gamma\sum\limits_{i=1}^k1_{B_\delta(z_{0,i})}\left(\mathcal{U}_{\gamma,\varepsilon,\mathbf{z_\varepsilon,a_\varepsilon}}-\mu_{\varepsilon,i}\right)_+^{\gamma-1}v.
\end{equation*}
By equation \eqref{1-5}, the error term $\omega_\varepsilon$ satisfies
\begin{equation}\label{3-5}
\mathbb{L}_{\gamma,\varepsilon}\omega_\varepsilon=H_\varepsilon+N_\varepsilon(\omega_\varepsilon),
\end{equation}
where
\begin{equation*}
H_\varepsilon=\sum\limits_{i=1}^k1_{B_\delta(z_{0,i})}\left(\mathcal{U}_{\gamma,\varepsilon,\mathbf{z_\varepsilon,a_\varepsilon}}-\mu_{\varepsilon,i}\right)_+^\gamma-\sum\limits_{i=1}^k1_{B_\delta(z_{0,i})}\left(W_{\gamma,\varepsilon,z_{\varepsilon,i},a_{\varepsilon,i}}-\frac{a_{\varepsilon,i}}{2\pi}\ln\frac{1}{\varepsilon}\right)_+^\gamma,
\end{equation*}
and
\begin{equation*}
	\begin{split}
    N_\varepsilon(v)&=\sum\limits_{i=1}^k1_{B_\delta(z_{0,i})}\left(\left(\mathcal{U}_{\gamma,\varepsilon,\mathbf{z_\varepsilon,a_\varepsilon}}+v-\mu_{\varepsilon,i}\right)_+^\gamma-\left(W_{\gamma,\varepsilon,z_{\varepsilon,i},a_{\varepsilon,i}}-\frac{a_{\varepsilon,i}}{2\pi}\ln\frac{1}{\varepsilon}\right)_+^\gamma\right)\\
    &\quad-\sum\limits_{i=1}^k1_{B_\delta(z_{0,i})}\gamma\left(W_{\gamma,\varepsilon,z_{\varepsilon,i},a_{\varepsilon,i}}-\frac{a_{\varepsilon,i}}{2\pi}\ln\frac{1}{\varepsilon}\right)_+^{\gamma-1}v.
    \end{split}
\end{equation*}
From estimates \eqref{3-2}, \eqref{3-4} and Proposition \eqref{prop2-2}, we see that 
\begin{equation}\label{3-6}
\mathbb{L}_{\gamma,\varepsilon}\omega_\varepsilon=H_\varepsilon+N_\varepsilon(\omega_\varepsilon)=0, \quad \text{in} \ \Omega\setminus\cup^k_{i=1}B_{Lr_{\varepsilon,i}(z_{\varepsilon,i})}.
\end{equation}

To derive a coercive estimate for the linear problem, we first discuss the limiting operator for $\mathbb L_{\gamma,\varepsilon}$ as $\varepsilon\to 0^+$, that is
\begin{equation*}
	\mathbb L_{\gamma,0} v=-\Delta v-\gamma (w_\gamma)_+^{\gamma-1}v, \quad \text{for} \ \gamma\in(0,1)\cup(1,\infty),
\end{equation*}
with $w_\gamma$ defined in \eqref{2-3}, and
\begin{equation*}
	\mathbb L_{1,0} v=-\Delta v-1_{B_\tau(0)}v, \quad \text{for} \ \gamma=1,
\end{equation*}
with $\tau>0$ defined before \eqref{2-4}. We have already known the following non-degenerate property for ground state $w_\gamma$ with $\gamma\in (0,\infty)$, see \cite{Cao1,DY,FW}.
\begin{proposition}\label{prop3-2}
	Let $w_\gamma$ with $\gamma\in(0,\infty)$ be the function defined in \eqref{2-3} and \eqref{2-4}. Suppose that $v\in L^\infty(\mathbb R^2)\cap C(\mathbb R^2)$ solves $\mathbb L_{\gamma,0} v=0$ in $\mathbb R^2$. Then
	$$v\in \mathrm{span}\,\left\{\frac{\partial w_\gamma}{\partial y_1},\frac{\partial w_\gamma}{\partial y_2}\right\}.$$
\end{proposition}

With this result in hand, we can prove the following lemma for $\gamma\in[1,\infty)$.
\begin{lemma}\label{lem3-3}
	Suppose that $\gamma\in [1,\infty)$ and $p\in (2,\infty)$. Then there is a $\varepsilon_0>0$ sufficiently small and some constants $\rho_0>0$, such that for any $\varepsilon\in (0,\varepsilon_0]$, one has
	\begin{equation*}
	\varepsilon^{-\frac{2}{p}}||\mathbb{L}_{\gamma,\varepsilon}\omega_\varepsilon||_{L^p(\cup^k_{i=1}B_{Lr_{\varepsilon,i}}(z_{\varepsilon,i}))}
	 \ge \rho_0||\omega_\varepsilon||_{L^\infty(\Omega)}.
	\end{equation*}
\end{lemma}
\begin{proof}
	Our strategy to derive the coercive estimate is to conduct a blow-up analysis and derive a contradiction. For the case $\gamma\in(1,\infty)$, we denote $\tilde{\omega}_{\varepsilon,i}(y)=\omega_\varepsilon(s_{\varepsilon,i}y+z_{\varepsilon,i})$ and
	$$\tilde{\mathbb{L}}_{\gamma,\varepsilon,i} v=-\Delta v-\gamma\varepsilon^{-2}\cdot\sum\limits_{i=1}^k s_{\varepsilon,i}^2\cdot\left(\mathcal{U}_{\gamma,\varepsilon,\mathbf{z_\varepsilon,a_\varepsilon}}(s_{\varepsilon,i}y+z_{\varepsilon,i})-\mu_{\varepsilon,i}\right)_+^{\gamma-1}v.$$
	Now we argue by contradiction. Suppose there are $\varepsilon_n\to 0^+$ such that $\varepsilon_n$ satisfies
	\begin{equation}\label{3-7}
	||\omega_{\varepsilon_n}||_{L^\infty(\Omega)}=1
	\end{equation}
	and
	\begin{equation*}
	\varepsilon_n^{-\frac{2}{p}}||\mathbb{L}_{\gamma,\varepsilon_n}\omega_{\varepsilon_n}||_{L^p(\cup^k_{i=1}B_{Ls_{\varepsilon_n,i}}(z_{\varepsilon_n,i}))}\le\frac{1}{n},
	\end{equation*}
    where we can use $B_{Ls_{\varepsilon_n,i}}(z_{\varepsilon_n,i})$ to substitute $B_{Ls_{\varepsilon_n,i}}(z_{\varepsilon_n,i})$ since it holds
    $|r_{i,\varepsilon}-s_{i,\varepsilon}|=o(\varepsilon)$ by \eqref{3-2}. Then after scaling, for $p\in(2,\infty)$ we have
	\begin{equation*}
	||\tilde{\omega}_{\varepsilon_n,i}||_{L^\infty(\Omega_{\gamma,\varepsilon_n,i})}=1,
	\end{equation*}
	and
	\begin{equation}\label{3-8}
	||\tilde{\mathbb{L}}_{\gamma,\varepsilon_n,i}\tilde{\omega}_{\varepsilon_n,i}||_{L^p(B_L(0))}\le\frac{c_i}{n},
	\end{equation}
	with $\Omega_{\gamma,\varepsilon_n,i}:=\{ s_{\varepsilon_n,i}y+z_{\varepsilon_n,i}\in\Omega\}$ and $c_i>0$ a constant independent of $\varepsilon$. By letting $g_n=\tilde{\mathbb{L}}_{\gamma,\varepsilon_n,i}\tilde{\omega}_{\varepsilon_n,j}$, we deduce that
	\begin{equation*}
	-\Delta\tilde{\omega}_{\varepsilon_n,i}=\gamma\varepsilon^{-2}\sum\limits_{i=1}^k s_{\varepsilon_n,i}^2\left(\mathcal{U}_{\gamma,\varepsilon_n,\mathbf{x}_{\varepsilon_n},\mathbf{a}_{\varepsilon_n}}(s_{\varepsilon_n,i}y+z_{\varepsilon_n,i})-\mu_{\varepsilon,i}\right)_+^{\gamma-1}\tilde{\omega}_{\varepsilon_n,i}+g_n.
	\end{equation*}
	where the right hand side is bounded in $L^p_{\text{loc}}(\mathbb{R}^2)$. According to the standard elliptic regularity theory, $\tilde{\omega}_{\varepsilon_n,i}$ is bounded in $W^{2,p}_{\text{loc}}(\mathbb{R}^2)$, and hence bounded in $C^{1}_{\text{loc}}(\mathbb{R}^2)$ for some $\alpha>0$ due to Sobolev embedding. As a result, we can assume that $\tilde{\omega}_{\varepsilon_n,i}$ converges uniformly in any compact set of $\mathbb{R}^2$ to $\omega^*\in L^\infty(\mathbb{R}^2)$. By \eqref{3-8}, the limiting function $\omega^*$ satisfies
	\begin{equation*}
	-\Delta \omega^*=\gamma (w_\gamma)_+^{\gamma-1}\omega^*, \quad \text{in} \ \mathbb{R}^2.
	\end{equation*}
	By Proposition \ref{prop3-2}, it must hold
	\begin{equation*}
	\omega^*=C_1\frac{\partial w_\gamma}{\partial y_1}+C_2\frac{\partial w_\gamma}{\partial y_2}.
	\end{equation*}
	On the other hand, since $p_{\varepsilon_n,i}$ is a local maximum point of $\psi_\varepsilon$, it holds $\nabla \omega_{\varepsilon_n}(p_{\varepsilon_n,i})=0$ and $\frac{p_{\varepsilon_n,i}-z_{\varepsilon_n,i}}{s_{\varepsilon_n,i}}\to 0$ by \eqref{3-2}. Hence we claim $\nabla\omega^*(0)=0$, which implies $C_1=C_2=0$ and $\omega^*\equiv 0$. As a result, we have actually proved $\omega_{\varepsilon_n}=o_\varepsilon(1)$ in $B_{Ls_{\varepsilon_n,i}}(z_{\varepsilon_n,i})$ for any large $L>0$. It follows from \eqref{3-6} that
	\begin{equation*}
	\mathbb{L}_{\gamma,\varepsilon_n}\omega_{\varepsilon_n}=0, \quad \text{in} \ \Omega\setminus\cup^k_{i=1}B_{Ls_{\varepsilon_n,i}(x_{\varepsilon_n,i})},
	\end{equation*}
	and
	\begin{equation*}
	-\Delta\omega_{\varepsilon_n}=0, \quad \text{in}  \ \Omega\setminus\cup^k_{i=1}B_{Ls_{\varepsilon_n,i}(x_{\varepsilon_n,i})}.
	\end{equation*}
	However, we know that $\omega_{\varepsilon_n}=0$ on $\partial\Omega$ and $\omega_{\varepsilon_n}=o_\varepsilon(1)$ on $\partial B_{Ls_{\varepsilon_n,i}}(x_{\varepsilon_n,i})$. By the maximum principle, it holds
	\begin{equation*}
	||\omega_{\varepsilon_n}||_{L^\infty(\Omega)}=o_\varepsilon(1), \ \ \ \text{as}  \ n\to\infty.
	\end{equation*} 
	Thus we have obtained a contradiction from \eqref{3-7}, and complete the proof for $\gamma\in(1,\infty)$. 
	
	For the case $\gamma=1$, the proof is similar, where we just need to substitute the scaling parameter $s_{\varepsilon,i}$ with $\varepsilon$. 
\end{proof}

When $\gamma\in (0,1)$, due to the singularity in linearized operator $\mathbb{L}_{\gamma,\varepsilon}$, we need the following estimate in $W^{-1,p}$ norm. Since the argument is similar to that in Lemma \ref{3-3}, we omit its proof (one can also use the argument for Proposition 3.3 in \cite{Cao4} for $\gamma=0$). 
\begin{lemma}\label{lem3-4}
	Suppose that $\gamma\in (0,1)$ and $p\in (2,\infty)$. Then there is a $\varepsilon_0>0$ sufficiently small and some constants $\rho_0>0$, such that for any $\varepsilon\in (0,\varepsilon_0]$, one has
	\begin{equation*}
		\varepsilon^{1-\frac{2}{p}}||\mathbb{L}_{\gamma,\varepsilon}\omega_\varepsilon||_{W^{-1,p}(\cup^k_{i=1}B_{Lr_{\varepsilon,i}}(z_{\varepsilon,i}))}
		\ge \rho_0||\omega_\varepsilon||_{L^\infty(\Omega)}.
	\end{equation*}
\end{lemma}

By now, we can give a more accurate estimate for $\omega_\varepsilon$ by the terms $\partial\mathcal W_{\gamma,i}$, $i=1,\cdots,k$ defined in the Appendix.
\begin{lemma}\label{lem3-5}
	For $\gamma\in(0,\infty)$, one has
	$$||\omega_\varepsilon||_{L^\infty(\Omega)}=O\left(\sum_{i=1}^k r_{\varepsilon,i}\partial\mathcal W_{\gamma,i}+\varepsilon^2\right),$$
	where $\partial\mathcal W_{\gamma,i}$ is defined in the Appendix. 
\end{lemma}
\begin{proof}
	Thanks to \eqref{3-5} and lemma \ref{lem3-3}, for $\gamma\in [1,\infty)$ and $p\in(2,\infty)$ we have
	\begin{equation}\label{3-9}
	||\omega_\varepsilon||_{L^\infty(\Omega)}
	\le C\varepsilon^{-\frac{2}{p}}||H_\varepsilon+N_\varepsilon(\omega_\varepsilon)||_{L^p(\cup^k_{i=1}B_{Lr_{\varepsilon,i}}(z_{\varepsilon,i}))}.
	\end{equation}
	and we only need to estimate $||H_\varepsilon+N_\varepsilon(\omega_\varepsilon)||_{L^p(\cup^k_{i=1}B_{Lr_{\varepsilon,i}}(z_{\varepsilon,i}))}$. So as to simplify our proof, in the following we will use $||\cdot||_p$ to denote the $L^p$ norm on $\cup^k_{i=1}B_{Lr_{\varepsilon,i}}(z_{\varepsilon,i})$, and use $||\cdot||_\infty$ to denote the $L^\infty$ norm $\Omega$. 
	
	Let us first consider the case $\gamma\in(1,\infty)$.
	From \eqref{A-1} in the Appendix and the H\"older inequality, we have
	\begin{equation*}
	\begin{split}
	||H_\varepsilon||_p&=\left\|\sum\limits_{i=1}^k1_{B_\delta(z_{0,i})}\left(\mathcal{U}_{\gamma,\varepsilon,\mathbf{z_\varepsilon,a_\varepsilon}}-\mu_{\varepsilon,i}\right)_+^\gamma-\sum\limits_{i=1}^k1_{B_\delta(z_{0,i})}\left(W_{\gamma,\varepsilon,z_{\varepsilon,i},a_{\varepsilon,i}}-\frac{a_{\varepsilon,i}}{2\pi}\ln\frac{1}{\varepsilon}\right)_+^\gamma\right\|_p\\
	&\le \varepsilon^{\frac{2}{p}}\cdot O\left(\sum_{i=1}^k s_{\varepsilon,i}\partial\mathcal W_{\gamma,i}+\varepsilon^2\right).
	\end{split}
	\end{equation*}
	While using Lemma \ref{A1} in the Appendix, we can derive
	\begin{equation*}
	\begin{split}
	||N_\varepsilon(\omega_\varepsilon)||_p&=\Bigg\|\sum\limits_{i=1}^k1_{B_\delta(z_{0,i})}\left(\left(\mathcal{U}_{\gamma,\varepsilon,\mathbf{z_\varepsilon,a_\varepsilon}}+\omega_\varepsilon-\mu_{\varepsilon,i}\right)_+^\gamma-\left(W_{\gamma,\varepsilon,z_{\varepsilon,i},a_{\varepsilon,i}}-\frac{a_{\varepsilon,i}}{2\pi}\ln\frac{1}{\varepsilon}\right)_+^\gamma\right)\\
	&\quad-\sum\limits_{i=1}^k1_{B_\delta(z_{0,i})}\gamma\left(W_{\gamma,\varepsilon,z_{\varepsilon,i},a_{\varepsilon,i}}-\frac{a_{\varepsilon,i}}{2\pi}\ln\frac{1}{\varepsilon}\right)_+^{\gamma-1}\omega_\varepsilon\Bigg\|_p\\
	&\le \varepsilon^{\frac{2}{p}}\cdot \|\omega_\varepsilon\|_\infty\cdot O\left(\|\omega_\varepsilon\|_\infty+\sum_{i=1}^k s_{\varepsilon,i}\partial\mathcal W_{\gamma,i}+\varepsilon^2\right)\\
	&=o(\varepsilon^{\frac{2}{p}})\cdot \|\omega_\varepsilon\|_\infty
	\end{split}
	\end{equation*}
	Then from \eqref{3-9}, we conclude
	$$||\omega_\varepsilon||_\infty\le o_\varepsilon(1)\cdot||\omega_\varepsilon||_\infty+O\left(\sum_{i=1}^k s_{\varepsilon,i}\partial\mathcal W_{\gamma,i}+\varepsilon^2\right),$$ 
	and hence we obtain the estimate for $\gamma\in(1,\infty)$.	
	
	For the case $\gamma=1$, we shall use the estimate \eqref{A-2} to deduce
	\begin{equation*}
		\begin{split}
			||H_\varepsilon||_p&=\left\|\sum\limits_{i=1}^k1_{B_\delta(z_{0,i})}\left(\mathcal{U}_{1,\varepsilon,\mathbf{z_\varepsilon,a_\varepsilon}}-\mu_{\varepsilon,i}\right)_+-\sum\limits_{i=1}^k1_{B_\delta(z_{0,i})}\left(W_{1,\varepsilon,z_{\varepsilon,i},a_{\varepsilon,i}}-\frac{a_{\varepsilon,i}}{2\pi}\ln\frac{1}{\varepsilon}\right)_+\right\|_p\\
			&\le \varepsilon^{\frac{2}{p}}\cdot O\left(\sum_{i=1}^k \tau\varepsilon\cdot \partial\mathcal W_{1,i}+\varepsilon^2\right).
		\end{split}
	\end{equation*}
    Since $N_\varepsilon(\omega_\varepsilon)$ takes the form
    \begin{equation*}
    	\begin{split}
    		N_\varepsilon(v)&=\sum\limits_{i=1}^k1_{B_\delta(z_{0,i})}\left(\left(\mathcal{U}_{1,\varepsilon,\mathbf{z_\varepsilon,a_\varepsilon}}+\omega_\varepsilon-\mu_{\varepsilon,i}\right)_+-\left(W_{1,\varepsilon,z_{\varepsilon,i},a_{\varepsilon,i}}-\frac{a_{\varepsilon,i}}{2\pi}\ln\frac{1}{\varepsilon}\right)_+\right)\\
    		&\quad-\sum\limits_{i=1}^k1_{B_\delta(z_{0,i})} 1_{B_{\tau\varepsilon}(z_{\varepsilon,i})}\omega_\varepsilon,
    	\end{split}
    \end{equation*}
    we will treat $N_\varepsilon$ in a way slightly different. We define
    $$S_{\varepsilon,i}^1:=\left(B_\delta(z_{0,i})\setminus B_{\tau\varepsilon}(z_{\varepsilon,i})\right) \, \Delta \, \{x\in B_\delta(z_{0,i})\,:\,\mathcal{U}_{1,\varepsilon,\mathbf{x_\varepsilon,a_\varepsilon}}(x)+\omega_\varepsilon<\mu_{\varepsilon,i}\},$$
    and
    $$S_{\varepsilon,i}^2:=B_{\tau\varepsilon}(z_{\varepsilon,i}) \, \Delta \, \{x\in B_\delta(z_{0,i})\,:\,\mathcal{U}_{1,\varepsilon,\mathbf{x_\varepsilon,a_\varepsilon}}(x)+\omega_\varepsilon<\mu_{\varepsilon,i}\},$$
    where $\Delta$ denotes the symmetric difference of two sets. By Lemma \ref{A2} in the appendix, we have
    \begin{equation*}
    	\cup_{i=1}^k|S_{\varepsilon,i}^1|+\cup_{i=1}^k|S_{\varepsilon,i}^2|=\varepsilon^2\cdot O\left(\sum_{i=1}^k \tau\varepsilon\cdot \partial\mathcal W_{1,i}+\varepsilon^2\right).
    \end{equation*}
    Hence we can derive
    \begin{equation*}
    	\begin{split}
    		||N_\varepsilon(\omega_\varepsilon)||_p&=\Bigg\|\sum\limits_{i=1}^k1_{B_\delta(z_{0,i})}\left(\left(\mathcal{U}_{1,\varepsilon,\mathbf{z_\varepsilon,a_\varepsilon}}+\omega_\varepsilon-\mu_{\varepsilon,i}\right)_+-\left(W_{1,\varepsilon,z_{\varepsilon,i},a_{\varepsilon,i}}-\frac{a_{\varepsilon,i}}{2\pi}\ln\frac{1}{\varepsilon}\right)_+\right)\\
    		&\quad-\sum\limits_{i=1}^k1_{B_\delta(z_{0,i})} 1_{B_{\tau\varepsilon}(z_{\varepsilon,i})}\omega_\varepsilon\Bigg\|_p\\
    		&\le C||\omega_\varepsilon||_\infty\cdot(\cup_{i=1}^k|S_{\varepsilon,i}^1|+\cup_{i=1}^k|S_{\varepsilon,i}^2|)^{\frac{1}{p}}\\
    		&=o(\varepsilon^{\frac{2}{p}})\cdot||\omega_\varepsilon||_\infty.
    	\end{split}
    \end{equation*}
    Using \eqref{3-9}, we obtain
    $$||\omega_\varepsilon||_\infty\le o_\varepsilon(1)\cdot||\omega_\varepsilon||_\infty+O\left(\sum_{i=1}^k \tau\varepsilon\cdot \partial\mathcal W_{1,i}+\varepsilon^2\right).$$ 
    So we obtain the desired result.
    
    Finally, we consider the case $\gamma\in(0,1)$. According to Lemma \ref{lem3-4}, it holds
    \begin{equation*}
    	||\omega_\varepsilon||_{L^\infty(\Omega)}
    	\le C\varepsilon^{1-\frac{2}{p}}||H_\varepsilon+N_\varepsilon(\omega_\varepsilon)||_{W^{-1,p}(\cup^k_{i=1}B_{Lr_{\varepsilon,i}}(z_{\varepsilon,i}))}
    \end{equation*}
    in this situation. It is easy to see that 
    $$||H_\varepsilon||_{W^{-1,p}(\cup^k_{i=1}B_{Lr_{\varepsilon,i}}(z_{\varepsilon,i}))}\le \varepsilon^{\frac{2}{p}-1}\cdot O\left(\sum_{i=1}^k \tau\varepsilon\cdot \partial\mathcal W_{\gamma,i}+\varepsilon^2\right).$$
    Moreover, following the strategy for the case $\gamma\in (1,\infty)$, we can derive
    \begin{equation*}
    	\begin{split}
    		&\quad||N_\varepsilon(\omega_\varepsilon)||_{W^{-1,p}(\cup^k_{i=1}B_{Lr_{\varepsilon,i}}(z_{\varepsilon,i}))}\\
    		&\le\left\|\sum\limits_{i=1}^k1_{B_\delta(z_{0,i})}\gamma\left(W_{\gamma,\varepsilon,z_{\varepsilon,i},a_{\varepsilon,i}}-\frac{a_{\varepsilon,i}}{2\pi}\ln\frac{1}{\varepsilon}\right)_+^{\gamma-1}\right\|_{W^{-1,p}(\cup^k_{i=1}B_{Lr_{\varepsilon,i}}(z_{\varepsilon,i}))}\\
    		&\quad\cdot\|\omega_\varepsilon\|_\infty\cdot O\left(\|\omega_\varepsilon\|_\infty+\sum_{i=1}^k s_{\varepsilon,i}\partial\mathcal W_{\gamma,i}+\varepsilon^2\right)\\
    		&= \varepsilon^{\frac{2}{p}-1}\cdot \|\omega_\varepsilon\|_\infty\cdot O\left(\|\omega_\varepsilon\|_\infty+\sum_{i=1}^k s_{\varepsilon,i}\partial\mathcal W_{\gamma,i}+\varepsilon^2\right)\\
    		&=o(\varepsilon^{\frac{2}{p}-1})\cdot \|\omega_\varepsilon\|_\infty.
    	\end{split}
    \end{equation*}
    Combining the estimates above, we have
    $$||\omega_\varepsilon||_\infty\le O\left(\sum_{i=1}^k s_{\varepsilon,i}\partial\mathcal W_{\gamma,i}+\varepsilon^2\right).$$ 
    Thus the proof is complete. 
\end{proof}

Thanks to Lemma \ref{lem3-5}, Lemma \ref{A1} and Lemma \ref{A2} in the Appendix can be rewritten as follows.
\begin{lemma}\label{lem3-6}
	The set 
	$$\Gamma_{\varepsilon,i}:=\{y \,:\, \tilde{\mathcal U}_{\gamma,\varepsilon,\mathbf{x_\varepsilon,a_\varepsilon},i}+\tilde\omega_{\varepsilon,i}=\mu_{\varepsilon,i}\}\cap B_L(0)$$
	is a closed curve in $\mathbb{R}^2$. 
	
	For $\gamma\in(0,1)\cup(1,\infty)$, $\Gamma_{\varepsilon,i}$ can be parameterized as
	\begin{equation*}
	\Gamma_{\varepsilon,i}(\theta)=(1+t_{\varepsilon,i})(\cos\theta,\sin\theta)=(\cos\theta,\sin\theta)+O\left(\sum_{i=1}^k s_{\varepsilon,i}\partial\mathcal W_{\gamma,i}+\varepsilon^2\right);
	\end{equation*}
    while for $\gamma=1$, $\Gamma_{\varepsilon,i}$ can be parameterized as
    \begin{equation*}
    	\Gamma_{\varepsilon,i}(\theta)=(\tau+t_{\varepsilon,i})(\cos\theta,\sin\theta)=(\cos\theta,\sin\theta)+O\left(\sum_{i=1}^k \tau\varepsilon\cdot \partial\mathcal W_{1,i}+\varepsilon^2\right),
    \end{equation*}
    with $\theta\in [0,2\pi)$.
\end{lemma}
 
With Lemma \ref{lem3-6} in hand, we can improve the estimates in Lemma \ref{lem2-1}.
\begin{lemma}\label{lem3-7}
	For any $x\in \Omega\setminus\cup^k_{i=1}B_{Lr_{\varepsilon,i}}(z_{\varepsilon,i})$, one has
	\begin{equation}\label{3-10}
	\psi_\varepsilon(x)=\sum\limits_{i=1}^k\kappa_i G(z_{\varepsilon,i},x)+O\left(\sum\limits_{i=1}^k\frac{\varepsilon \sup_{\theta\in[0,2\pi)}|t_{\varepsilon,i}(\theta)|}{|x-z_{\varepsilon,i}|} +\sum\limits_{i=1}^k\frac{ \varepsilon^2}{|x-z_{\varepsilon,i}|^2}\right),
	\end{equation}
	and 
	\begin{equation}\label{3-11}
	\frac{\partial \psi_\varepsilon(x)}{\partial x_h}=\sum\limits_{i=1}^k\kappa_j \frac{\partial G(z_{\varepsilon,i},x)}{\partial x_h}+O\left(\sum\limits_{i=1}^k\frac{\varepsilon \sup_{\theta\in[0,2\pi)}|t_{\varepsilon,i}(\theta)|}{|x-z_{\varepsilon,i}|} +\sum\limits_{i=1}^k\frac{ \varepsilon^2}{|x-z_{\varepsilon,i}|^2}\right)
	\end{equation}	
    with $h=1,2$.
\end{lemma}
\begin{proof}
	For any $x\in \Omega\setminus\cup^k_{i=1}B_{Lr_{\varepsilon,i}}(z_{\varepsilon,i})$, using a same expansion in the proof of Lemma \ref{lem2-1}, we can deduce
	\begin{equation*}
	\begin{split}
	\psi_\varepsilon(x)&=\frac{1}{\varepsilon^2}\sum\limits_{i=1}^k1_{B_\delta(z_{0,i})}\int_{\Omega}(\psi_\varepsilon(x)-\mu_{\varepsilon,i})_+^\gamma G(y,x)dy\\
	&=\sum\limits_{i=1}^k\kappa_i G(z_{\varepsilon,i},x)+\frac{1}{\varepsilon^2}\sum\limits_{i=1}^k1_{B_\delta(z_{0,i})}\int_{\Omega}(\psi_\varepsilon(x)-\mu_{\varepsilon,i})_+^\gamma(G(y,x)-G(z_{\varepsilon,i},x))dy\\
	&=\sum\limits_{i=1}^k\kappa_i G(z_{\varepsilon,i},x)+\frac{1}{\varepsilon^2}\sum\limits_{i=1}^k\int_{D_{\varepsilon,i}}(\psi_\varepsilon(x)-\mu_{\varepsilon,i})_+^\gamma\langle\nabla G(z_{\varepsilon,i},x),y-z_{\varepsilon,i}\rangle dy\\
	& \ \ \ +O(\sum\limits_{i=1}^k\frac{\varepsilon^2}{|x-z_{\varepsilon,i}|}).\\
	\end{split}
	\end{equation*}
	By Lemma \ref{lem3-6} and odd symmetry, for $i=1,\cdots,k$ we have
	\begin{equation*}
	\begin{split}
	&\int_{D_{\varepsilon,i}}(\psi_\varepsilon(x)-\mu_{\varepsilon,i})_+^\gamma\langle\nabla G(z_{\varepsilon,i},x),y-z_{\varepsilon,i}\rangle dy\\
	&=\int_{D_{\varepsilon,i}\setminus B_{Lr_{\varepsilon,i}}(z_{\varepsilon,i})}(\psi_\varepsilon(x)-\mu_{\varepsilon,i})_+^\gamma\langle\nabla G(z_{\varepsilon,i},x),y-z_{\varepsilon,i}\rangle dy\\
	& \quad +\int_{B_{Lr_{\varepsilon,i}}(z_{\varepsilon,i})}(\mathcal{U}_{\gamma,\varepsilon,\mathbf{x_\varepsilon,a_\varepsilon}}+\omega_\varepsilon-\mu_{\varepsilon,i})_+^\gamma\langle\nabla G(z_{\varepsilon,i},x),y-z_{\varepsilon,i}\rangle dy\\
	&=O\left(\frac{\varepsilon}{|x-z_{\varepsilon,i}|}\cdot|D_{\varepsilon,i}\setminus B_{Lr_{\varepsilon,i}}(x_{\varepsilon,j})|\right)+O_\varepsilon(1)\cdot\int_{B_{Lr_{\varepsilon,i}}(z_{\varepsilon,i})}\omega_\varepsilon\langle\nabla G(z_{\varepsilon,i},x),y-z_{\varepsilon,i}\rangle dy\\
	&=O\left(\sum\limits_{j=1}^k\frac{\varepsilon^3\sup_{\theta\in[0,2\pi)}|t_{\varepsilon,j}(\theta)|}{|x-z_{\varepsilon,j}|}+\frac{ \varepsilon^2}{|x-z_{\varepsilon,i}|^4}\right).
	\end{split}
	\end{equation*}	
    Hence \eqref{3-10} is obvious and we can obtain \eqref{3-11} by a similar strategy.
\end{proof}

Now we can bootstrap to get our desired estimates.
\begin{proposition}\label{prop3-8}
	Suppose $i=1,\cdots,k$. The following estimates hold
	\begin{equation}\label{3-12}
	||\omega_\varepsilon||_{L^\infty(\Omega)}=O(\varepsilon^2),
	\end{equation}
	\begin{equation}\label{3-13}
		||t_{\varepsilon,i}||_{L^{\infty}([0,2\pi))}=O(\varepsilon^2),
	\end{equation}
	\begin{equation}\label{3-14}
		|\mathbf z_\varepsilon-\mathbf z_0|=O(\varepsilon^2).
	\end{equation}
	For $\gamma\in(0,1)\cup(1,\infty)$, one has
	\begin{equation*}
		s_{\varepsilon,i}=\varepsilon\cdot\left(\frac{2\pi\phi'_\gamma(1)}{\kappa_i}\right)^{\frac{\gamma-1}{2}}+O(\varepsilon^3)=s_{0,i}+O(\varepsilon^3),
	\end{equation*}
	\begin{equation*}
		\mu_{\varepsilon,i}=\frac{\kappa_i}{2\pi} \ln \frac{1}{\varepsilon}+\frac{\gamma-1}{4\pi}\cdot\kappa_i \ln \frac{2\pi\phi_\gamma'(1)}{\kappa_i}-\kappa_iH(z_{0,i},z_{0,i})+\sum_{j\neq i}\kappa_j G(z_{0,i},z_{0,j})+O(\varepsilon^2|\ln\varepsilon|).
	\end{equation*}
    For $\gamma=1$, one has
    \begin{equation*}
    	\Lambda_{\varepsilon,i}=\frac{\kappa_i}{2\pi\tau\cdot\phi'_1(\tau)}+O(\varepsilon^2)=\Lambda_{0,i}+O(\varepsilon^2),
    \end{equation*}
    \begin{equation*}
    	\mu_{\varepsilon,i}=\frac{\kappa_i}{2\pi} \ln \frac{1}{\varepsilon}+\frac{\kappa_i}{2\pi}\ln\tau-\kappa_iH(z_{0,i},z_{0,i})+\sum_{j\neq i}\kappa_j G(z_{0,i},z_{0,j})+O(\varepsilon^2|\ln\varepsilon|).
    \end{equation*}
\end{proposition}
\begin{proof}
	We will use an argument of bootstrap to obtain the estimates. At the first stage, by the definition of $W_{\gamma,i}$ in the Appendix, we have $\partial\mathcal W_{\gamma,i}=O_\varepsilon(1)$ for $i=1,\cdots,k$. Using Lemma \ref{lem3-5} and Lemma \ref{lem3-6}, we have
	\begin{equation*}
		||\omega_\varepsilon||_{L^\infty(\Omega)}=O(\varepsilon) \quad \text{and} \quad ||t_{\varepsilon,i}||_{L^{\infty}([0,2\pi))}=O(\varepsilon), \quad \mathrm{for} \ i=1,\cdots,k.
	\end{equation*}
    Then, by the Pohozaev identity, we have
    $$\int_{B_1(0)}\phi_\gamma^\gamma=2\pi|\phi_\gamma'(1)|, \quad \text{for} \ \gamma\in(0,1)\cup(1,\infty),$$
    and 
    $$\int_{B_\tau(0)}\phi_1=2\pi\tau|\phi_1'(\tau)|, \quad \text{for} \ \gamma=1.$$
    Recall condition $(\mathbf{B})$ on local circulation $\kappa_i$. For $\gamma\in(0,1)\cup(1,\infty)$, one has
    \begin{equation*}
	   s_{\varepsilon,i}=\varepsilon\cdot\left(\frac{2\pi\phi'_\gamma(1)}{\kappa_i+O(\varepsilon)}\right)^{\frac{\gamma-1}{2}}+O(\varepsilon^2)=s_{0,i}+O(\varepsilon^2),
	\end{equation*}
    and
	 \begin{equation*}
		a_{\varepsilon,i}\cdot\frac{|\ln\varepsilon|}{|\ln s_{\varepsilon,i}|}=\left(\frac{\varepsilon}{s_{\varepsilon,i}}\right)^{\frac{2}{\gamma-1}}2\pi\phi'_\gamma(1)=\kappa_i+O(\varepsilon)
	\end{equation*}
    with $i=1,\cdots,k$. Hence for $\gamma\in(0,1)\cup(1,\infty)$, we obtain
    \begin{equation}\label{3-15}
    \partial\mathcal W_{\gamma,i}=(1+O(\varepsilon))\cdot\partial_{\mathbf x_i}\mathcal W_k(\mathbf z_\varepsilon), \quad \text{for} \ i=1,\cdots,k.
    \end{equation}
    Similarly, for $\gamma=1$ one has
    \begin{equation*}
    	\Lambda_{\varepsilon,i}=\frac{\kappa_i+O(\varepsilon)}{2\pi\tau\cdot\phi'_1(\tau)}=\Lambda_{0,i}+O(\varepsilon),
    \end{equation*}
    and
    \begin{equation*}
    	a_{\varepsilon,i}\cdot\frac{|\ln\varepsilon|}{|\ln \tau\varepsilon|}=2\pi\tau\Lambda_{\varepsilon,i}\cdot\phi'_1(\tau)=\kappa_i+O(\varepsilon)
    \end{equation*}
    with $i=1,\cdots,k$. Thus it holds
    \begin{equation}\label{3-16}
    	\partial\mathcal W_{1,i}=(1+O(\varepsilon))\cdot\partial_{\mathbf x_i}\mathcal W_k(\mathbf z_\varepsilon), \quad \text{for} \ i=1,\cdots,k.
    \end{equation}
	
	On the other hand, using Lemma \ref{lem3-7} and the Pohozaev identity \eqref{1-7}, we can proceed as the proof of Theorem \ref{thm1} to obtain
	$$\partial_{\mathbf x_i}\mathcal W_k(\mathbf z_\varepsilon)=O(\varepsilon^2), \quad \text{for} \ i=1,\cdots,k,$$
	which yields \eqref{3-14} by the nondegeneracy of $\mathbf z_0=(z_{0,1},\cdots,z_{0,k})$. Combining this estimate with \eqref{3-15} and \eqref{3-16}, we have
	$$\partial\mathcal W_{\gamma,i}=O(\varepsilon^2), \quad \text{for} \ i=1,\cdots,k \ \text{and} \ \gamma\in(0,\infty).$$
	Thus we have improved the estimate for $\partial\mathcal W_{\gamma,i}$ from $O_\varepsilon(1)$ to $O(\varepsilon^2)$. Then we can repeat our first step to conclude \eqref{3-12} and \eqref{3-13}.
	
	For $\gamma\in(0,1)\cup(1,\infty)$, it is easy to see that
	\begin{equation*}
		s_{\varepsilon,i}=\varepsilon\cdot\left(\frac{2\pi\phi'_\gamma(1)}{\kappa_i}\right)^{\frac{\gamma-1}{2}}+O(\varepsilon^3)=s_{0,i}+O(\varepsilon^3)
	\end{equation*}
	for $i=1,\cdots,k$. Using equation \eqref{3-1}, we have
	\begin{equation*}
	\mu_{\varepsilon,i}=\frac{\kappa_i}{2\pi} \ln \frac{1}{\varepsilon}+\frac{\gamma-1}{4\pi}\cdot\kappa_i \ln \frac{2\pi\phi_\gamma'(1)}{\kappa_i}-\kappa_iH(z_{0,i},z_{0,i})+\sum_{j\neq i}\kappa_j G(z_{0,i},z_{0,j})+O(\varepsilon^2|\ln\varepsilon|).
	\end{equation*}
	For the case $\gamma=1$, we can calculate in parallel to derive 
	\begin{equation*}
	 	\Lambda_{\varepsilon,i}=\frac{\kappa_i}{2\pi\tau\cdot\phi'_1(\tau)}+O(\varepsilon^2)=\Lambda_{0,i}+O(\varepsilon^2),
	\end{equation*}
    and
    \begin{equation*}
	\mu_{\varepsilon,i}=\frac{\kappa_i}{2\pi} \ln \frac{1}{\varepsilon}+\frac{\kappa_i}{2\pi}\ln\tau-\kappa_iH(z_{0,i},z_{0,i})+\sum_{j\neq i}\kappa_j G(z_{0,i},z_{0,j})+O(\varepsilon^2|\ln\varepsilon|).
	\end{equation*}
	Hence we have verified the last part of this Proposition and finished our proof.
\end{proof}

\section{The proof of uniqueness}

In this section, we will prove the local uniqueness of vortices stated in Theorem \ref{thm2}. Our strategy is arguing by contradiction. Suppose that \eqref{1-5} together with conditions $(\mathbf{A})$ and $(\mathbf{B})$ has two different solutions $\psi_\varepsilon^{(1)}$ and $\psi_\varepsilon^{(2)}$, which blow up at $\mathbf z_0$. Let $x_{\varepsilon,i}^{(m)}$, $\mu_{\varepsilon,i}^{(m)}$ and so on be the parameters or functions appearing in $\psi_\varepsilon^{(i)}$, and
\begin{equation*}
\zeta_\varepsilon(x)=\frac{\psi_\varepsilon^{(1)}(x)-\psi_\varepsilon^{(2)}(x)}{||\psi_\varepsilon^{(1)}-\psi_\varepsilon^{(2)}||_\infty}.
\end{equation*}
be the normalized difference, where we simplify the norm $\|\cdot\|_{L^\infty(\Omega)}$ by $\|\cdot\|_\infty$. Then $\zeta_\varepsilon$ must satisfy $||\xi_\varepsilon||_\infty$ and the equation
\begin{equation*}
	\begin{cases}
		-\Delta\zeta_\varepsilon(x)=f_\varepsilon(x), &x\in\Omega,\\
		\zeta_\varepsilon(x)=0, &x\in\partial \Omega, \\
	\end{cases}
\end{equation*}
where
\begin{equation*}
f_\varepsilon(x)=\frac{1}{||\psi_\varepsilon^{(1)}-\psi_\varepsilon^{(2)}||_\infty}\sum\limits_{i=1}^k\left(\mathbf 1_{B_\delta(x_{0,i})}(\psi_\varepsilon^{(1)}-\mu_{\varepsilon,i}^{(1)})_+^\gamma-\mathbf1_{B_\delta(z_{0,i})}(\psi_\varepsilon^{(2)}-\mu_{\varepsilon,i}^{(2)})_+^\gamma\right).
\end{equation*}
Since it holds
$$-\Delta\zeta_\varepsilon=0,\quad \mathrm{in} \ \Omega\setminus\cup^k_{i=1}B_{Lr_{\varepsilon,i}}(z_{0,i})$$
by estimates given in Proposition \ref{prop3-8}. To obtain a contradiction and verify the uniqueness, we only need to prove $$\zeta_\varepsilon=o_\varepsilon(1), \quad \mathrm{in} \ \cup^k_{i=1}B_{Lr_{\varepsilon,i}}(z_{0,i})$$
and apply the maximum principle.

In the following lemma, we will first prove a limiting result for $\gamma\in [1,\infty)$, where we use subscript $(1)$ to denote the parameters in $\psi_\varepsilon^{(1)}$.
\begin{lemma}\label{lem4-1}
	Suppose $p\in(2,\infty)$ and $i=1,\cdots,k$. Then as $\varepsilon\to 0^+$, for $\gamma\in (1,\infty)$, one has  $\|s_{\varepsilon,i}^2f_\varepsilon(s_{\varepsilon,i}y+x_{\varepsilon,i}^{(1)})\|_{L^p(B_L(0))}\le C$, 
	\begin{equation}\label{4-1}
		s_{\varepsilon,i}^2f_\varepsilon(s_{\varepsilon,i}y+z_{\varepsilon,i}^{(1)})\to \gamma(w_\gamma)_+^{\gamma-1}\zeta_i, \quad \mathrm{in} \ L^p_{\mathrm{loc}}(B_L(0));
	\end{equation}
    while for $\gamma=1$, one has $||\varepsilon^2f_\varepsilon(\varepsilon y+z_{\varepsilon,i}^{(1)})||_{L^p(B_L(0))}\le C$,
    \begin{equation}\label{4-2}
    	\varepsilon^2f_\varepsilon(\varepsilon y+z_{\varepsilon,i}^{(1)})\to 1_{B_\gamma(0)}\zeta_i, \quad \mathrm{in} \ L^p_{\mathrm{loc}}(B_L(0)),
    \end{equation}
	where $w_\gamma$ is defined in \eqref{2-3} \eqref{2-4}, and $\zeta_i$ is the limiting function satisfying 
	\begin{equation*}
	\zeta_\varepsilon(s_{\varepsilon,i}y+z_{\varepsilon,i}^{(1)})\to\zeta_i=b_{1,i}\frac{\partial w_\gamma}{\partial y_1}+b_{2,i}\frac{\partial w_\gamma}{\partial y_2}, \quad \mathrm{in} \ C_{\mathrm{loc}}^{1,\alpha}(B_L(0)),
	\end{equation*}
	as $\varepsilon\to 0^+$ with $b_{1,i}$ and $b_{2,i}$ two constants.
\end{lemma}
\begin{proof}
	According to the estimates in Proposition \ref{prop3-8}, we have 
	\begin{equation*}
	\begin{split}
	f_\varepsilon(x)&=\frac{1}{\varepsilon^2||\psi_\varepsilon^{(1)}-\psi_\varepsilon^{(2)}||_\infty}\sum\limits_{i=1}^k1_{B_\delta(z_{0,i})}\left((\mathcal{U}_{\gamma,\varepsilon,\mathbf{z_\varepsilon,a_\varepsilon}}^{(1)}+\omega_\varepsilon^{(1)}-\mu_{\varepsilon,i}^{(1)})_+^\gamma-(\mathcal{U}_{\gamma,\varepsilon,\mathbf{z_\varepsilon,a_\varepsilon}}^{(2)}+\omega_\varepsilon^{(2)}-\mu_{\varepsilon,i}^{(2)})_+^\gamma\right)\\
	&=\frac{1}{\varepsilon^2||\psi_\varepsilon^{(1)}-\psi_\varepsilon^{(2)}||_\infty}\bigg(\sum\limits_{i=1}^k1_{B_\delta(z_{0,i})}\big(W_{\gamma,\varepsilon,z_{\varepsilon,i},a_{\varepsilon,i}}^{(1)}(x)-\frac{a_{\varepsilon,i}^{(1)}}{2\pi}\ln\frac{1}{\varepsilon}+\omega_\varepsilon^{(1)}+O(\varepsilon^2)\big)_+^\gamma\\
	& \quad -\sum\limits_{i=1}^k1_{B_\delta(z_{0,i})}\big(W_{\gamma,\varepsilon,z_{\varepsilon,i},a_{\varepsilon,i}}^{(2)}(x)-\frac{a_{\varepsilon,i}^{(2)}}{2\pi}\ln\frac{1}{\varepsilon}+\omega_\varepsilon^{(2)}+O(\varepsilon^2)\big)_+^\gamma\bigg)\\
	&=\frac{\gamma}{\varepsilon^2||\psi_\varepsilon^{(1)}-\psi_\varepsilon^{(2)}||_\infty}\sum\limits_{i=1}^k1_{B_\delta(z_{0,i})}\big(W^{(1)}_{\gamma,\varepsilon,z_{\varepsilon,i},a_{\varepsilon,i}}-\frac{a^{(1)}_{\varepsilon,i}}{2\pi}\ln\frac{1}{\varepsilon}+O(\varepsilon^2)\big)_+^{\gamma-1}\\
	&\quad\cdot(\psi_\varepsilon^{(1)}-\psi_\varepsilon^{(2)}-(\mu^{(1)}_{\varepsilon,i}-\mu^{(2)}_{\varepsilon,i})).
	\end{split}
	\end{equation*}
    In view of the circulation constraint ($\mathbf{B}$), we deduce that
    \begin{equation*}
    	\begin{split}
    	0&=\int_{B_\delta(z_{0,i})}\left((\psi_\varepsilon^{(1)}-\mu_{\varepsilon,i}^{(1)})_+^\gamma-(\psi_\varepsilon^{(2)}-\mu_{\varepsilon,i}^{(2)})_+^\gamma\right)dx\\
    	&=\int_{B_\delta(z_{0,i})}\big(W^{(1)}_{\gamma,\varepsilon,z_{\varepsilon,i},a_{\varepsilon,i}}-\frac{a^{(1)}_{\varepsilon,i}}{2\pi}\ln\frac{1}{\varepsilon}+O(\varepsilon^2)\big)_+^{\gamma-1}(\psi_\varepsilon^{(1)}-\psi_\varepsilon^{(2)}-(\mu^{(1)}_{\varepsilon,i}-\mu^{(2)}_{\varepsilon,i}))dx,
    	\end{split}
    \end{equation*}
    which implies 
    $$\frac{\mu_{\varepsilon,i}^{(1)}-\mu_{\varepsilon,i}^{(2)}}{||\psi_\varepsilon^{(1)}-\psi_\varepsilon^{(2)}||_\infty}\le C$$
    for $i=1,\cdots,i$. Then using the fact $\|\zeta_\varepsilon\|_\infty=1$ and expansion for $f_\varepsilon(x)$ above, we see that
	$$||s_{\varepsilon,i}^2f_\varepsilon(s_{\varepsilon,i}y+z_{\varepsilon,i}^{(1)})||_{L^p(B_L(0))}\le C, \quad \mathrm{for} \ p\in(2,\infty).$$
    By $L^p$ estimate for elliptic operator, $\zeta_\varepsilon(s_{\varepsilon,i}y+z_{\varepsilon,i}^{(1)})$ is bounded in $W_{\text{loc}}^{2,p}(\mathbb{R}^2)$. In view of Sobolev embedding, as $\varepsilon\to0^+$, we may assume $\zeta_\varepsilon(s_{\varepsilon,i}y+z_{\varepsilon,i}^{(1)})\to\zeta_i$ in $C_{\text{loc}}^{1,\alpha}(\mathbb{R}^2)$. Hence it is obvious that
    \begin{equation*}
    \begin{split}
    s_{\varepsilon,i}^2f_\varepsilon(s_{\varepsilon,i}y+z_{\varepsilon,i}^{(1)})&=\gamma (w_\gamma)_+^{\gamma-1}\zeta_\varepsilon(s_{\varepsilon,i}y+z_{\varepsilon,i}^{(1)})+O(\varepsilon^2)\\
    &\to \gamma (w_\gamma)_+^{\gamma-1}\zeta_i, \quad \text{in} \ L^p_{\text{loc}}(\mathbb{R}^2).
    \end{split}
    \end{equation*}
	As a result, we claim that \eqref{4-1} holds and $\zeta_i$ satisfies
	$$-\Delta\zeta_i=\gamma (w_\gamma)_+^{\gamma-1}\zeta_i, \quad \text{in} \ \mathbb{R}^2.$$
	This gives
	$$\zeta_i=b_{1,i}\frac{\partial w_\gamma}{\partial y_1}+b_{2,i}\frac{\partial w_\gamma}{\partial y_2}.$$
	
	When $\gamma=1$, we can proceed as the former case to derive \eqref{4-2}. Hence the proof is complete.
\end{proof}

If $\gamma\in(0,1)$, there will be singularities in the linearized operator. However, by introducing the $W^{-1,p}$ norm and using the regularity theory for elliptic operator, we can derive a similar limiting result just as Lemma \ref{lem4-1} (see also Lemma 4.1 in \cite{Cao4} for the case $\gamma=0$). Recalling that $\Omega_{\gamma,\varepsilon,i}:=\{y\in\mathbb{R}^2:s_{\varepsilon,i}y+z_{\varepsilon,i}^{(1)}\in\Omega\}$ denotes the domain after scaling, we have following lemma concerning this situation.
\begin{lemma}\label{lem4-2}
	Suppose $p\in(2,\infty)$ and $i=1,\cdots,k$. Then as $\varepsilon\to 0^+$, for $\gamma\in (0,1)$, one has $\|s_{\varepsilon,i}^2f_\varepsilon(s_{\varepsilon,i}y+x_{\varepsilon,i}^{(1)})\|_{W^{-1,p}(\Omega_{\gamma,\varepsilon,i})}\le C$, 
	\begin{equation*}
		\int_{\Omega_{\gamma,\varepsilon,i}}s_{\varepsilon,i}^2f_\varepsilon(s_{\varepsilon,i} y+z_{\varepsilon,i}^{(1)})\varphi dy\to \int\gamma(w_\gamma)_+^{\gamma-1}\zeta_i \varphi dy, \quad \forall \, \phi\in C^\infty_0(\mathbb{R}^2).
	\end{equation*}
    where $w_\gamma$ is defined in \eqref{2-3}, and $\zeta_i$ is the limiting function satisfying 
	\begin{equation*}
	\zeta_\varepsilon(\varepsilon y+z_{\varepsilon,i}^{(1)})\to\zeta_i=b_{1,i}\frac{\partial w_\gamma}{\partial y_1}+b_{2,i}\frac{\partial w_\gamma}{\partial y_2}, \quad \mathrm{in} \ C_{\mathrm{loc}}^\alpha(\mathbb{R}^2),
	\end{equation*}
	as $\varepsilon\to 0$ with $b_{1,i}$ and $b_{2,i}$ two constants.
\end{lemma}

To verify $\|\zeta_\varepsilon\|_\infty=o_\varepsilon(1)$ and derive a contradiction, we need to prove $b_{1,i}=b_{2,i}=0$ for $i=1,\cdots,k.$. To achieve this goal, we will use a local Pohozaev identity technique. The following lemma gives an essential estimate for $\zeta_\varepsilon$ in $\Omega\setminus \cup_{i=1}^kB_{\delta}(z_{\varepsilon,i}^{(1)})$.
\begin{lemma}\label{lem5-3}
    It holds
	\begin{equation}\label{5-8}
	\zeta_\varepsilon(x)=\sum\limits_{i=1}^k\sum\limits_{l=1}^2B_{l,i,\varepsilon}\partial_hG(z_{\varepsilon,i}^{(1)},x)+O(\varepsilon^2), \quad \mathrm{in} \ C^1(\Omega\setminus \cup_{i=1}^k B_{\delta}(z_{\varepsilon,i}^{(1)})),
	\end{equation}
	with $\partial_l G(y,x)= \partial G(y,x)/\partial y_l$, and
	\begin{equation*}
	B_{l,i,\varepsilon}=
	\int_{B_\delta(z_{\varepsilon,i}^{(1)})}\left(y_l-z_{\varepsilon,i,l}^{(1)}\right)f_{\varepsilon}(y)dy.
	\end{equation*}
\end{lemma}
\begin{proof}
	We have
	\begin{equation*}
	\begin{split}
	\zeta_\varepsilon(x)=&\int_{\Omega}G(y,x)f_\varepsilon(y)dy=\sum\limits_{i=1}^k\int_{B_{Lr_{\varepsilon,i}}(z_{\varepsilon,i}^{(1)})}G(y,x)f_\varepsilon(y)dy\\
	=&\sum\limits_{i=1}^kA_{\varepsilon,i}G(z_{\varepsilon,i}^{(1)},x)+\sum\limits_{i=1}^k\sum\limits_{h=1}^2B_{h,i,\varepsilon}\partial_hG(z_{\varepsilon,i}^{(1)},x)\\
	&+\sum\limits_{i=1}^k\int_{B_{Lr_{\varepsilon,i}}(z_{\varepsilon,i}^{(1)})} \left(G(y,x)-G(z_{\varepsilon,i}^{(1)},x)-\langle\nabla G(z_{\varepsilon,i}^{(1)},x),y-z_{\varepsilon,i}^{(1)} \rangle \right) f_\varepsilon(y)dy,
	\end{split}
	\end{equation*}
	where
	$$A_{\varepsilon,i}:=\int_{B_{Lr_{\varepsilon,i}}(z_{\varepsilon,i}^{(1)})}f_\varepsilon(y)dy=0.$$
	On the other hand, for $\gamma\in(0,1)\cup(1,\infty)$ it holds
	$$\int_{B_{Lr_{\varepsilon,i}}(z_{\varepsilon,i}^{(1)})}|f_\varepsilon(y)|dy=\int_{B_{L}(0)}s_{\varepsilon,i}^2|f_\varepsilon(s_{\varepsilon,i}y+z_{\varepsilon,i}^{(1)})|dy
	\le C,$$
	while for $\gamma=1$ it holds
	$$\int_{B_{Lr_{\varepsilon,i}}(z_{\varepsilon,i}^{(1)})}|f_\varepsilon(y)|dy=\int_{B_{L}(0)}\varepsilon^2|f_\varepsilon(\varepsilon y+z_{\varepsilon,i}^{(1)})|dy
	\le C,$$
	Hence we deduce by the estimates in Proposition \ref{prop3-8} that
	\begin{equation*}
	\begin{split}
	&\left|\int_{B_{Lr_{\varepsilon,i}}(z_{\varepsilon,i}^{(1)})}\left(G(y,x)-G(z_{\varepsilon,i}^{(1)},x)-\langle\nabla G(z_{\varepsilon,i}^{(1)},x),y-z_{\varepsilon,i}^{(1)} \rangle \right) f_\varepsilon(y)dy\right|\\
	\leq& C \varepsilon^2\int_{B_{Lr_{\varepsilon,i}}(z_{\varepsilon,i}^{(1)})}|f_\varepsilon(y)|dy
	=O(\varepsilon^2),
	\end{split}
	\end{equation*}
	and thus \eqref{5-8} follows. 
\end{proof}

Having these preparations done, we can go ahead for the uniqueness of solutions to \eqref{1-5} satisfying assumptions $(\mathbf{A})$ and $(\mathbf{B})$.

{\bf Proof of Theorem \ref{thm2}:} By applying \eqref{1-7} on $\psi_\varepsilon^{(1)}$, $\psi_\varepsilon^{(2)}$ separately and make the difference, we have following local Pohozaev identity for $\zeta_\varepsilon$: For $i=1,\cdots,k$ and $0<\delta_0<\delta$, it holds
\begin{equation*}
-\int_{\partial B_{\delta_0}(z_{\varepsilon,i}^{(1)})}\frac{\partial\zeta_\varepsilon}{\partial\nu}\frac{\psi_\varepsilon^{(1)}}{\partial x_h}-\int_{\partial B_{\delta_0}(z_{\varepsilon,i}^{(1)})}\frac{\partial \psi_\varepsilon^{(2)}}{\partial\nu}\frac{\partial\zeta_\varepsilon}{\partial x_h}+\frac{1}{2}\int_{\partial B_{\delta_0}(x_{\varepsilon,i}^{(1)})}\langle\nabla(\psi_\varepsilon^{(1)}+\psi_\varepsilon^{(2)}),\nabla\zeta_\varepsilon\rangle\nu_h=0
\end{equation*}
with $h=1,2$. By Lemma \ref{lem5-3}, we obtain 
\begin{equation*}
\begin{split}
&-\int_{\partial B_{\delta_0}(z_{\varepsilon,i}^{(1)})}\frac{\partial\zeta_\varepsilon}{\partial\nu}\sum\limits_{j=1}^k\kappa_jD_{x_h}G(x_{\varepsilon,j}^{(1)},x)-\int_{\partial B_{\delta_0}(z_{\varepsilon,i}^{(1)})}\langle\sum\limits_{j=1}^k\kappa_j\nabla_y G(z_{\varepsilon,j}^{(1)},x),\nu\rangle\frac{\partial\zeta_\varepsilon}{\partial\nu}\\
&+\frac{1}{2}\int_{\partial B_{\delta_0}(z_{\varepsilon,i}^{(1)})}\sum\limits_{j=1}^k\kappa_l\langle D_xG(x_{\varepsilon,j}^{(1)},x),\nabla_y\zeta_\varepsilon\rangle\nu_h=O(\varepsilon^2).
\end{split}
\end{equation*}
According to Proposition \ref{prop3-8}, the above identity can be rewritten as
\begin{equation}\label{4-4}
\begin{split}
&-\int_{\partial B_{\delta_0}(z_{\varepsilon,i}^{(1)})}\sum\limits_{j=1}^k\sum\limits_{m=1}^k\sum\limits_{l=1}^2\kappa_jB_{l,m,\varepsilon}\langle D_x\partial_{y_l}G(z_{\varepsilon,m}^{(1)},x),\nu\rangle D_{x_h}G(z_{\varepsilon,i}^{(1)},x)\\
&-\int_{\partial B_{\delta_0}(z_{\varepsilon,i}^{(1)})}\sum\limits_{j=1}^k\sum\limits_{m=1}^k\sum\limits_{l=1}^2\kappa_jB_{l,m,\varepsilon}\langle \nabla_y G(z_{\varepsilon,j}^{(1)},x),\nu\rangle D_{x_h}\partial_{y_l} G(z_{\varepsilon,m}^{(1)},x)\\
&+\int_{\partial B_{\delta_0}(z_{\varepsilon,i}^{(1)})}\sum\limits_{j=1}^k\sum\limits_{m=1}^k\sum\limits_{l=1}^2\kappa_jB_{l,m,\varepsilon}\langle \nabla_y G(z_{\varepsilon,j}^{(1)},x),D_{x_h}\partial_{y_l} G(z_{\varepsilon,m}^{(1)},x)\rangle\nu_h=O(\varepsilon^2).
\end{split}
\end{equation}
We define the following quadratic form:
\begin{equation*}
Q(u,v)=-\int_{\partial B_{\delta_0}(z_{\varepsilon,i}^{(1)})}\frac{\partial v}{\partial \nu}\frac{\partial u}{\partial x_h}
-\int_{\partial B_{\delta_0}(z_{\varepsilon,i}^{(1)})}\frac{\partial u}{\partial \nu}\frac{\partial v}{\partial x_h}
+\int_{\partial B_{\delta_0}(z_{\varepsilon,i}^{(1)})}\langle\nabla u, \nabla v\rangle \nu_h.
\end{equation*}
Note that if $u$ and $v$ are harmonic in $B_\delta(z_{\varepsilon,j}^{(1)})\setminus \{z_{\varepsilon,j}^{(1)}\}$, then $Q(u,v)$ is independent of $\delta_0\in (0,\delta]$. It is easy to check that if $j\neq i$ and $m\neq i$,
\begin{equation*}
Q(G(z_{\varepsilon,j}^{(1)},x),\partial_{y_l}G(z_{\varepsilon,m}^{(1)},x))=0.
\end{equation*}
So from \eqref{4-4}, we obtain
\begin{equation}\label{4-5}
\begin{split}
&\sum\limits_{l=1}^2\kappa_iQ(G(z_{\varepsilon,i}^{(1)},x),\partial_{y_l}G(z_{\varepsilon,i}^{(1)},x))B_{l,i,\varepsilon}+\sum\limits_{j\neq i}\sum\limits_{l=1}^2\kappa_jQ(G(z_{\varepsilon,j}^{(1)},x),\partial_{y_l}G(z_{\varepsilon,i}^{(1)},x))B_{l,i,\varepsilon}\\
&+\sum\limits_{m\neq i}\sum\limits_{l=1}^2\kappa_iQ(G(z_{\varepsilon,j}^{(1)},x),\partial_{y_l}G(z_{\varepsilon,m}^{(1)},x))B_{l,m,\varepsilon}=O(\varepsilon^2).
\end{split}
\end{equation}
By denoting $B_{l,m,\varepsilon}=\kappa_m \tilde{B}_{l,m,\varepsilon}$,  $\mathbf{B}_\varepsilon=(\tilde{B}_{1,1,\varepsilon},\tilde{B}_{2,1,\varepsilon},\cdots,\tilde{B}_{1,k,\varepsilon},\tilde{B}_{2,k,\varepsilon})$, and using the local Pohozaev identities for Green's function (Appendix 6.2 in \cite{CPYb}), we have
\begin{equation*}
\begin{split}
&\sum\limits_{l=1}^2\kappa_iQ(G(z_{\varepsilon,i}^{(1)},x),\partial_{y_l}G(z_{\varepsilon,i}^{(1)},x))B_{l,i,\varepsilon}+\sum\limits_{j\neq i}\sum\limits_{l=1}^2\kappa_jQ(G(z_{\varepsilon,j}^{(1)},x),\partial_{y_l}G(z_{\varepsilon,i}^{(1)},x))B_{l,i,\varepsilon}\\
&+\sum\limits_{m\neq i}\sum\limits_{l=1}^2\kappa_iQ(G(z_{\varepsilon,j}^{(1)},x),\partial_{y_l}G(z_{\varepsilon,m}^{(1)},x))B_{l,m,\varepsilon}=-\frac{1}{2}D^2_x\mathcal{W}_k(z_{\varepsilon,1}^{(1)},\cdots,z_{\varepsilon,k}^{(1)})\cdot\mathbf{B}_\varepsilon.
\end{split}
\end{equation*}
Then from \eqref{4-5}, we find
\begin{equation*}
D^2_x\mathcal{W}_k(z_{\varepsilon,1}^{(1)},\cdots,z_{\varepsilon,k}^{(1)})\cdot\mathbf{B}_\varepsilon=O(\varepsilon^2).
\end{equation*}
This relation together the non-degeneracy of the critical point $\mathbf{z_0}$, implies
$$B_{l,i,\varepsilon}=O(\varepsilon^2), \quad \ \text{for} \ l=1,2 \ \text{and} \ i=1,\cdots,k.$$
However, it holds
$$B_{l,i,\varepsilon}=\int_{B_{Lr_{\varepsilon,i}}(z_{\varepsilon,i}^{(1)})}(y_l-z_{\varepsilon,i,l}^{(1)})f_{\varepsilon}(y)dy=O_\varepsilon\left(\int_{B_1(0)}(b_{1,i}\frac{\partial w_\gamma}{\partial y_1}+b_{2,i}\frac{\partial w_\gamma}{\partial y_2})\cdot y_ldy\right).$$
Thus, $b_{1,i}=b_{2,i}=0$. So we have proved $|\zeta_\varepsilon|=o_\varepsilon(1)$ in $B_{Lr_{\varepsilon,i}}(z_{\varepsilon,i}^{(1)})$. On the other hand, since 
$$-\Delta\zeta_\varepsilon=0, \quad \text{in} \ \Omega\setminus\cup_{i=1}^kB_{Lr_{\varepsilon,i}}(z_{\varepsilon,i}^{(1)})$$
with zero boundary condition on $\partial\Omega$. Using the maximum principle, we conclude $\|\zeta_\varepsilon\|_\infty=o_\varepsilon(1)$, which is a contradiction to $||\zeta_\varepsilon||_\infty=1$. Combining the local uniqueness with asymptotic estimates in Proposition \ref{prop3-8}, we have completed the proof of Theorem \ref{thm2}. \qed

\section{The nonlinear stability}

The topic of this section originates from Lord Kelvin's work \cite{K}, which told us that vortices of greatest energy relative to an isovortical surface constitute stable flows. This important observation was then developed by Arnol'd in \cite{Ar1,Ar2,Ar3}. In \cite{Bu2}, Burton reinterpreted the original idea of Kelvin, and obtained a criterion for nonlinear stability of steady vortices in a bounded domain: the vorticity $\omega$ corresponding to isolated energy maximizers is nonlinearly stable. In the following, we will show that the local uniqueness result obtained in Section 4 will lead to the isolation condition, and verify the nonlinear stability for a class of steady vortices in $\Omega$.

To illustrate our idea, we first give a variational characteristic for $\omega_\varepsilon=-\Delta\psi_\varepsilon$, where $\psi_\varepsilon$ is a solution to \eqref{1-5} satisfying conditions $(\mathbf{A})$, $(\mathbf{B})$, and $(\mathbf{C})$. Recall that the kinetic energy of a flow in $\Omega$ is defined by
$$E[\omega]=\frac{1}{2}\int_\Omega |\mathbf u|^2 dx=\frac{1}{2}\int_\Omega\omega\mathcal G\omega(x) dx.$$
For the purpose of dealing with the nonlinear term in the first equation of \eqref{1-5}, we also introduce a penalty term
$$\mathcal P_\varepsilon[\omega]=\frac{\gamma}{1+\gamma}\cdot\varepsilon^2\int_\Omega\omega^{1+\frac{1}{\gamma}}(x)dx.$$
Denote the admissible class $\mathcal A_\varepsilon$ of vorticity as
$$\mathcal A_\varepsilon=\left\{\omega\in L^\infty(\Omega) \,:\,  0\le\omega\le\frac{M}{\varepsilon^2}, \ \int_\Omega \omega(x)dx=\kappa_1\right\},$$
where $M>0$ is a sufficiently large constant independent of $\varepsilon$. We consider the following variational problem
\begin{equation}\label{5-1}
	\mathcal E_\varepsilon=\sup_{\omega\in\mathcal A_\varepsilon}\left(E[\omega]-\mathcal P_\varepsilon[\omega]\right),
\end{equation}
and let $\mathcal S_\varepsilon\subset \mathcal A_\varepsilon$ be the set of maximizers of \eqref{5-1}.

The following proposition is a corollary of Theorem 1.2 in \cite{CWZ}, where a relationship of solutions to \eqref{1-5} and variational problem \eqref{5-1} is established.
\begin{proposition}\label{prop5-1}
	If $\varepsilon\in (0,\varepsilon_0]$ with $\varepsilon_0>0$ sufficiently small, then $\mathcal S_\varepsilon\neq\emptyset$ and each maximizer $\omega_\varepsilon\in \mathcal S_\varepsilon$ gives a steady solution to \eqref{1-2}, which takes the form
	\begin{equation*}
		\omega_\varepsilon=\frac{1}{\varepsilon^2}(\mathcal G\omega_\varepsilon-\mu_\varepsilon)_+^\gamma, \quad \int_\Omega\omega_\varepsilon(x)dx=\kappa_1,
	\end{equation*}  
	where $\mu_\varepsilon$ is the flux constant satisfying $\mu_\varepsilon=-\frac{\kappa}{2\pi}\ln\varepsilon+O_\varepsilon(1)$. Moreover, as $\varepsilon\to 0$, the support of $\omega_\varepsilon$ shrinks to a minimum $z_*$ of the Robin function $\mathcal R(x)$ in the sense that 
	$$\sup_{x\in \mathrm{supp}\,\omega_\varepsilon}|x-z_*|\to 0, \quad \mathrm{diam}\, \mathrm{supp}\, \omega_\varepsilon\le R_0\varepsilon$$ with $R_0$ some positive constant.
\end{proposition}

Let $u$ be a non-negative Lebesgue integrable function in $\Omega$, we denote by $\mathcal F(u)$ the set of (equimeasurable) rearrangements of $u$ in $\Omega$ defined by
\begin{equation*}
	\mathcal F(u)=\Big\{ v\in L^1(\Omega):v\ge 0~ \text{and}~ \text{meas}\, \{x: v(x)>\tau\}=\text{meas}\, \{x: u(x)>\tau\}, \forall\, \tau>0  \Big\}.
\end{equation*}
Kelvin's notation of `isovortical surface' is somehow consistent with the rearrangement class $\mathcal F(\omega)$ of some prescribed vorticity $\omega$. To show the nonlinear stability stated in Theorem \ref{thm3}, we need the following Burton's criterion (Theorem 1 in \cite{Bu2}). 
\begin{proposition}\label{prop5-2}
	Let $\Omega\subset\mathbb{R}^2$ be a bounded simply connected domain of class $C^{2,\alpha}$ for some $\alpha\in(0,1)$, and $4/3<p<\infty$. Let $\bar\omega\in L^p(\Omega)$ be the isolated maximizer of kinetic energy $E$ over $\mathcal R(\bar\omega)$, namely, there exists $\xi_0>0$ such that for any $\omega\in \mathcal F(\bar\omega)$, $0<\|\omega-\bar\omega\|_{L^p(\Omega)}<\xi_0$, it holds $E(\bar\omega)>E(\omega)$. Then $\bar\omega$ is a steady solution of \eqref{1-2}, and nonlinearly stable in the following sense: for any $\eta>0$, there exists $\xi>0$ such that if $\omega(\cdot,t)\in L^{\infty}_{\mathrm{loc}}(\mathbb R, L^p(\Omega))$ with $4/3<p<\infty$ is an energy-conserving solution of \eqref{1-2}, and $\|\omega(0,\cdot)-\omega_\varepsilon\|_{L^p(\Omega)}<\xi$, then $\|\omega(t,\cdot)-\omega_\varepsilon\|_{L^p(\Omega)}<\eta$ for all $t\in\mathbb R$.
\end{proposition}

Now we can prove our nonlinear stability result.

{\bf Proof of Theorem \ref{thm3}:} Since $\omega_\varepsilon$ is a maximizer to \eqref{5-1} by Proposition \ref{prop5-1}, we can use the fact that $\mathcal F(\omega_\varepsilon)\subset \mathcal A_\varepsilon$, and deduce that $\omega_\varepsilon$ is a maximizer of $E-\mathcal P_\varepsilon$ over $\mathcal F(\omega_\varepsilon)$. However, $\mathcal P_\varepsilon$ is a constant on $\mathcal F(\omega_\varepsilon)$. Then according to Proposition \ref{prop5-2}, we only need to show that $\omega_\varepsilon$ is an isolated maximum of $E$ over $\mathcal A_\varepsilon$.

For $4/3<p<\infty$, let $\hat\omega_\varepsilon$ be another maximizer of $E$ over $\mathcal F(\omega_\varepsilon)$ satisfying $\|\hat\omega_\varepsilon-\omega_\varepsilon\|_{L^p(\Omega)}<\xi$. Then $\hat\omega_\varepsilon$ is also a maximizer to \eqref{5-1} by fine property of $\mathcal F(\omega_\varepsilon)$. In view of Proposition \ref{prop5-1}, $\hat\psi_\varepsilon=\mathcal G\hat\omega_\varepsilon$ is a solution to 
\begin{equation*}
	\begin{cases}
		-\varepsilon^2\Delta \hat\psi_\varepsilon=\mathbf1_{B_\delta(\hat z_{0,1})}(\hat\psi_\varepsilon-\hat\mu_\varepsilon)_+^\gamma, &x\in\Omega,\\
		\hat\psi_\varepsilon=0, &x\in\partial \Omega, \\
		\int_\Omega \frac{1}{\varepsilon^2}(\hat\psi_\varepsilon-\hat\mu_\varepsilon)_+^\gamma dx=\kappa_1,\\
		\mathrm{supp}\, (\hat\psi_\varepsilon-\hat\mu_\varepsilon)_+\subset B_{o_\varepsilon(1)}(\hat z_{0,1}),\\
		\mathrm{diam}\, \mathrm{supp}\,(\hat\psi_\varepsilon-\hat\mu_\varepsilon)_+\le R_0\varepsilon,
	\end{cases}
\end{equation*}
where $\hat z_{0,1}$ is a minimum point of the Robin function $\mathcal R(x)$ that may be different from $z_{0,1}$ in condition $(\mathbf{C})$. If $\hat z_{0,1}=z_{0,1}$, it holds $\hat\omega_\varepsilon=\omega_\varepsilon$ by Theorem \ref{thm2}. If $\hat z_{0,1}\neq z_{0,1}$, then by the fact that $z_{0,1}$ is a non-degenerate minimum point of $\mathcal R(x)$ and condition $(\mathbf{A})$, we see $\text{supp}\, \hat\omega_\varepsilon\cup \text{supp}\, \omega_\varepsilon=\emptyset$ when $\varepsilon$ is sufficiently small. Hence $\|\hat\omega_\varepsilon-\omega_\varepsilon\|_{L^p(\Omega)}=2\sigma$ with $\sigma$ the $L^p$ norm of elements in $\mathcal R(\omega_\varepsilon)$, which is a contradiction. So we have verified all assumptions of Proposition \ref{prop5-2} in our situation, and complete the proof of Theorem \ref{thm3}.\qed 

To finish the discussion of this section, we outline the proof of nonlinear stability in a more general situation, where $k>1$ is a finite number, and $\mathbf{z}_0=(z_{0,1},\cdots,z_{0,k})\in \Omega^k$ is a non-degenerate minimum point of $\mathcal{W}_k(\mathbf x)$. We introduce the admissible class $\mathcal A_\varepsilon^k$ as
$$\mathcal A_\varepsilon^k=\left\{\omega\in L^\infty(\Omega) \,:\,  \omega=\sum_{i=1}^k\omega_i, \ \text{supp}\, \omega_i\subset B_\delta(z_{0,i}), \ 0\le\omega_i\le\frac{M}{\varepsilon^2}, \ \int_\Omega \omega_i(x)dx=\kappa_i\right\}.$$
By Theorem 1.5 in \cite{CWZ}, each maximizer $\omega_\varepsilon$ of \eqref{5-1} in $\mathcal A_\varepsilon^k$ is a solution to \eqref{1-2}. Moreover, $\psi_\varepsilon=\mathcal G\omega_\varepsilon$ satisfies \eqref{1-5} together with conditions $(\mathbf{A})$ and $(\mathbf{B})$. In view of Theorem \ref{thm2}, $\omega_\varepsilon$ is unique. Then we can use a same argument as in the proof of Theorem \ref{thm3} to verify the isolation of each $\omega_{\varepsilon,i}$ in $B_\delta(z_{0,i})$, and obtain the nonlinear stability of $\omega_\varepsilon$ by Burton's criterion.

\bigskip

{\bf Acknowledgements:} The authors are grateful to the anonymous referee for helpful comments that have improved the exposition and clarity of this manuscript. The authors also thank Prof. Shusen Yan and Dr. Guolin Qin for their kind suggestions. This work was supported by NNSF of China Grant 11831009.

\bigskip

\appendix

\section{Estimates for the vorticity set}

In this appendix, we will give some estimates and statements for the vorticity set that have been used in the previous sections. For $\gamma\in (0,1)\cup(1,\infty)$, we denote $\tilde u_i=u(s_{\varepsilon,i}y+z_{\varepsilon,i})$. In view of \eqref{3-2}, for $i=1,\cdots,k$, and $y$ satisfying $|y|\le L$, we have
\begin{equation}\label{A-1}
	\begin{split}
		&\quad\tilde{\mathcal U}_{\gamma,\varepsilon,\mathbf{x_\varepsilon,a_\varepsilon},i}(y)-\mu_{\varepsilon,i}\\
		&=\tilde W_{\gamma,\varepsilon,i,i}(y)-\frac{a_{\varepsilon,i}}{2\pi}\ln\frac{1}{\varepsilon}-\frac{a_{\varepsilon,i}|\ln\varepsilon|}{|\ln s_{\varepsilon,i}|}\cdot \langle \nabla_xH(z_{\varepsilon,i},z_{\varepsilon,i}),s_{\varepsilon,i}y\rangle\\
		&\quad+\sum\limits_{j\neq i}\frac{a_{\varepsilon,j}|\ln\varepsilon|}{|\ln s_{\varepsilon,j}|}\cdot\langle \nabla_xG(z_{\varepsilon,i},z_{\varepsilon,j}),s_{\varepsilon,i}y\rangle+O(\varepsilon^2),
	\end{split}
\end{equation}
For $i=1,\cdots,k$, let
$$\mathcal N_{\gamma,\varepsilon,i}:=\left(\frac{\varepsilon}{s_{\varepsilon,i}}\right)^{\frac{2}{\gamma-1}}\phi'_\gamma(1)$$
be the absolute value of gradient for $\tilde W_{\gamma,\varepsilon,i,i}(y)$ on $|y|=1$, and 
$$\partial\mathcal W_{\gamma,i}:=\sum\limits_{j\neq i}\frac{a_{\varepsilon,j}|\ln\varepsilon|}{|\ln s_{\varepsilon,j}|}\cdot \nabla_xG(z_{\varepsilon,i},z_{\varepsilon,j})-\frac{a_{\varepsilon,i}|\ln\varepsilon|}{|\ln s_{\varepsilon,i}|}\cdot \nabla_xH(z_{\varepsilon,i},z_{\varepsilon,i}).$$
We have following lemma concerning the shape of $\partial D_{\varepsilon,i}$.
\begin{lemma}\label{A1}
	Suppose $\gamma\in(0,1)\cup(1,\infty)$, and
	\begin{equation}\label{A-2}
	||\omega_\varepsilon||_{L^\infty(\Omega)}+\varepsilon||\nabla\omega_\varepsilon||_{L^\infty(\Omega)}=o_\varepsilon(1).
	\end{equation}
	Then the set 
	$$\Gamma_{\varepsilon,i}:=\{y \,:\, \tilde{\mathcal U}_{\gamma,\varepsilon,\mathbf{x_\varepsilon,a_\varepsilon},i}+\tilde\omega_{\varepsilon,i}=\mu_{\varepsilon,i}\}\cap B_L(0)$$
	is a closed curve in $\mathbb{R}^2$, which can be written as 
	\begin{equation}\label{A-3}
	   \begin{split}
	   \Gamma_{\varepsilon,i}(\theta)&=(1+t_{\varepsilon,i})(\cos\theta,\sin\theta)\\
	   &=\left(1+\frac{1+o_\varepsilon(1)}{\mathcal N_{\gamma,\varepsilon,i}}\cdot\tilde\omega_{\varepsilon,i}(\cos\theta,\sin\theta)\right)(\cos\theta,\sin\theta)\\
	   &\quad+\frac{1+o_\varepsilon(1)}{\mathcal N_{\gamma,\varepsilon,i}}\cdot s_{\varepsilon,i}\partial\mathcal W_{\gamma,i}\cdot(\cos\theta,\sin\theta)\\
	   &\quad+O(\varepsilon^2), \quad \theta\in [0,2\pi).
	   \end{split}
	\end{equation}
	Moreover, it holds
	\begin{equation*}
	(\tilde{\mathcal{U}}_{\gamma,\varepsilon,\mathbf{x_\varepsilon,a_\varepsilon},i}+\tilde\omega_{\varepsilon,i})((1+t_{\varepsilon,i})(\cos\theta,\sin\theta))-\mu_{\varepsilon,i}\left\{
	\begin{array}{lll}
    >0, & \mathrm{if} \ t<t_{\varepsilon,i}(\theta),\\
	<0, & \mathrm{if} \ t>t_{\varepsilon,i}(\theta).
	\end{array}
	\right.
	\end{equation*}	
\end{lemma}
\begin{proof}
	Note that
	\begin{equation*}
	\tilde W_{\gamma,\varepsilon,i,i}(y)-\frac{a_{\varepsilon,i}}{2\pi}\ln\frac{1}{\varepsilon}=\left\{
	\begin{array}{lll}
	 (\frac{\varepsilon}{s_{\varepsilon,i}})^{\frac{2}{\gamma-1}}\phi_\gamma(|y|), & |y|\le 1,\\
	 a_{\varepsilon,i}\frac{|\ln\varepsilon|}{|\ln s_{\varepsilon,i}|}\ln\frac{1}{|y|}, & |y|\ge 1.
	\end{array}
	\right.
	\end{equation*}
	According to \eqref{A-1} and \eqref{A-2}, we find for arbitrarily small $\delta_1>0$, if $|y|<1-\delta_1$, then
	\begin{equation*}
	\begin{split}
	\tilde{\mathcal{U}}_{\gamma,\varepsilon,\mathbf{x_\varepsilon,a_\varepsilon,}j}(y)+\tilde\omega_{\varepsilon,i}(y)-\mu_{\varepsilon,i}&=(\frac{\varepsilon}{s_{\varepsilon,i}})^{\frac{2}{\gamma-1}}\phi_\gamma(|y|)+\tilde\omega_{\varepsilon,i}(y)+O(\varepsilon)\\
	&>(\frac{\varepsilon}{s_{\varepsilon,i}})^{\frac{2}{\gamma-1}}\phi_\gamma(|1-\delta_1|)+o_\varepsilon(1)>0.
	\end{split}
	\end{equation*}
	On the other hand, for arbitrarily small $\delta_2>0$, if $1+\delta_2<|y|<L$, then
	\begin{equation*}
	\begin{split}
	\tilde{\mathcal{U}}_{\gamma,\varepsilon,\mathbf{x_\varepsilon,a_\varepsilon,}j}(y)+\tilde\omega_{\varepsilon,i}(y)-\mu_{\varepsilon,i}&=a_{\varepsilon,i}\frac{|\ln\varepsilon|}{|\ln s_{\varepsilon,i}|}\ln\frac{1}{|y|}+\tilde\omega_{\varepsilon,i}(y)+O(\varepsilon)\\
	&<a_{\varepsilon,i}\frac{|\ln\varepsilon|}{|\ln s_{\varepsilon,i}|}\ln\frac{1}{|1+\delta_2|}+o_\varepsilon(1)<0.
	\end{split}
	\end{equation*}
	Thus we have verified that for any $(\cos\theta,\sin\theta)$, there exist a $t_{\varepsilon,i}(\theta)$, such that $|t_{\varepsilon,i}(\theta)|=o_\varepsilon(1)$, and
	$$(1+t_{\varepsilon,i})(\cos\theta,\sin\theta)\in\Gamma_{\varepsilon,i}(\theta)$$
	
	From Proposition \ref{prop2-2} and Lemma \ref{lem3-1} we see that
	\begin{equation*}
		\frac{\partial (\tilde{\mathcal{U}}_{\gamma,\varepsilon,\mathbf{x_\varepsilon,a_\varepsilon},i}+\tilde\omega_{\varepsilon,i})((1+t)(\cos\theta,\sin\theta))}{\partial t}\Bigg|_{t=t_{\varepsilon,i}}=-\mathcal N_{\gamma,\varepsilon,i}+o_\varepsilon(1)<0.
	\end{equation*}
	As a result, $t_{\varepsilon,i}$ is unique by the implicit function theorem and $ \Gamma_{\varepsilon,i}$ is a continuous closed curve in $\mathbb{R}^2$. Moreover, using \eqref{A-1} and the implicit function theorem again, we can derive \eqref{A-3} by direct calculation. Hence the proof of Lemma \ref{A1} is complete.
\end{proof}

Next, we consider the case $\gamma=1$, and denote $\tilde u_i=u(\varepsilon y+z_{\varepsilon,i})$. Using \eqref{3-4}, for $i=1,\cdots,k$, and $y$ satisfying $|y|\le L$, it holds
\begin{equation}\label{A-4}
	\begin{split}
		&\quad\tilde{\mathcal U}_{1,\varepsilon,\mathbf{x_\varepsilon,a_\varepsilon},i}(y)-\mu_{\varepsilon,i}\\
		&=\tilde W_{1,\varepsilon,i,i}(y)-\frac{a_{\varepsilon,i}}{2\pi}\ln\frac{1}{\varepsilon}-\frac{a_{\varepsilon,i}|\ln\varepsilon|}{|\ln \tau\varepsilon|}\cdot \langle \nabla_xH(z_{\varepsilon,i},z_{\varepsilon,i}),s_{\varepsilon,i}y\rangle\\
		&\quad+\sum\limits_{j\neq i}\frac{a_{\varepsilon,j}|\ln\varepsilon|}{|\ln\tau\varepsilon|}\cdot\langle \nabla_xG(z_{\varepsilon,i},z_{\varepsilon,j}),s_{\varepsilon,i}y\rangle+O(\varepsilon^2),
	\end{split}
\end{equation}
For $i=1,\cdots,k$, let
$$\mathcal N_{1,\varepsilon,i}:=\tau\Lambda_{\varepsilon,i}\cdot\phi'_1(\tau)$$
be the gradient constant on $|y|=\tau$, and 
$$\partial\mathcal W_{1,i}:=\sum\limits_{j\neq i}\frac{a_{\varepsilon,j}|\ln\varepsilon|}{|\ln\tau\varepsilon|}\cdot \nabla_xG(z_{\varepsilon,i},z_{\varepsilon,j})-\frac{a_{\varepsilon,i}|\ln\varepsilon|}{|\ln \tau\varepsilon|}\cdot \nabla_xH(z_{\varepsilon,i},z_{\varepsilon,i}).$$
We have following lemma as an analogy of Lemma \ref{A1}.  Since the estimates share a similar strategy, we omit its proof here.
\begin{lemma}\label{A2}
	Suppose $\gamma=1$, and
	\begin{equation*}
		||\omega_\varepsilon||_{L^\infty(\Omega)}+\varepsilon||\nabla\omega_\varepsilon||_{L^\infty(\Omega)}=o_\varepsilon(1).
	\end{equation*}
	Then the set 
	$$\Gamma_{\varepsilon,i}:=\{y \,:\, \tilde{\mathcal U}_{1,\varepsilon,\mathbf{x_\varepsilon,a_\varepsilon},i}+\tilde\omega_{\varepsilon,i}=\mu_{\varepsilon,i}\}\cap B_L(0)$$
	is a closed curve in $\mathbb{R}^2$, which can be written as 
	\begin{equation}\label{A-5}
		\begin{split}
			\Gamma_{\varepsilon,i}(\theta)&=(\tau+t_{\varepsilon,i})(\cos\theta,\sin\theta)\\
			&=\left(\tau+\frac{1+o_\varepsilon(1)}{\mathcal N_{1,\varepsilon,i}}\cdot\tilde\omega_{\varepsilon,i}(\cos\theta,\sin\theta)\right)(\cos\theta,\sin\theta)\\
			&\quad+\frac{1+o_\varepsilon(1)}{\mathcal N_{1,\varepsilon,i}}\cdot \tau\varepsilon\cdot \partial\mathcal W_{1,i}\cdot(\cos\theta,\sin\theta)\\
			&\quad+O(\varepsilon^2), \quad \theta\in [0,2\pi).
		\end{split}
	\end{equation}
	Moreover, it holds
	\begin{equation*}
		(\tilde{\mathcal{U}}_{1,\varepsilon,\mathbf{x_\varepsilon,a_\varepsilon},i}+\tilde\omega_{\varepsilon,i})((\tau+t_{\varepsilon,i})(\cos\theta,\sin\theta))-\mu_{\varepsilon,i}\left\{
		\begin{array}{lll}
			>0, & \mathrm{if} \ t<t_{\varepsilon,i}(\theta),\\
			<0, & \mathrm{if} \ t>t_{\varepsilon,i}(\theta).
		\end{array}
		\right.
	\end{equation*}	
\end{lemma}

\phantom{s}
 \thispagestyle{empty}

\end{document}